\PassOptionsToPackage{numbers,sort&compress}{natbib}  
\documentclass[preprint,12pt]{elsarticle}
\usepackage[margin=2.5cm]{geometry}
\usepackage{enumitem}
\newlist{inlinelist}{enumerate*}{1}
\setlist*[inlinelist,1]{%
  label=(\roman*),
}
\usepackage{dirtytalk}
\usepackage{soul}
\usepackage{hhline}
\usepackage{subfig}
\usepackage{multirow}
\usepackage{xcolor}
\usepackage{graphicx}
\usepackage{upgreek}
\usepackage{amssymb}
\usepackage{amsfonts,amsthm,bm,amsmath} 
\usepackage{appendix}
\usepackage{times}
\usepackage{caption}
\usepackage{subfig}
\newcommand{\psubref}[1]{\protect\subref{#1}}
\usepackage{multirow}
\usepackage{xcolor}
\usepackage[linesnumbered,ruled,vlined]{algorithm2e}

\SetCommentSty{mycommfont}

\SetKwInput{KwInput}{Input}                
\SetKwInput{KwOutput}{Output}              

\usepackage[colorlinks]{hyperref} 
\hypersetup{ 
    colorlinks=true,       
    linkcolor=red,          
    citecolor=blue,        
    filecolor=magenta,      
    urlcolor=cyan           
}

\newcommand{\fref}[1]{Fig.~\ref{#1}}

\newcommand{\eref}[1]{Eq.~(\ref{#1})}

\newcommand{\sref}[1]{Section~\ref{#1}}
\newcommand{\tref}[1]{Table~\ref{#1}}
\setcitestyle{square}

\journal{Structural and Multidisciplinary Optimization}
\begin{document}

\begin{frontmatter}

\title{LatticeOPT: A heuristic topology optimization framework for thin-walled, 2D extruded lattices}
\author[]{Junyan He$^1$}
\author[]{Shashank Kushwaha$^1$}
\author[]{Diab Abueidda$^2$}
\author[]{Iwona Jasiuk$^1$\corref{mycorrespondingauthor}}
\address{$^1$ Department of Mechanical Science and Engineering, University of Illinois at Urbana-Champaign, Champaign, IL, USA \\
$^2$ National Center for Supercomputing Applications, University of Illinois at Urbana-Champaign, Champaign, IL, USA}
\cortext[mycorrespondingauthor]{Corresponding author}
\ead{ijasiuk@illinois.edu}
\begin{abstract}

This paper introduces a heuristic topology optimization framework for thin-walled, 2D extruded lattice structures subject to complex high-speed loading. The proposed framework optimizes the wall thickness distribution in the lattice cross section through different thickness update schemes, inspired by the idea of equalization of absorbed energy density across all lattice walls. The proposed framework is ubiquitous and can be used in explicit dynamic simulations, which is the primary numerical method used in crashworthiness studies. No information on the material tangent stiffness matrix is required, and complex material behaviors and complex loading conditions can be handled. Three numerical examples are presented to demonstrate framework capabilities: (1) Optimization of a long, slender column under axial compression to maximize specific energy absorption, (2) Optimization of a lattice-filled sandwich panel under off-center blast loading to minimize material damage, (3) Generation of a periodic lattice core design under blast loading. The results show that the framework can effectively increase specific energy absorption or minimize material damage with as few as 25 finite element simulations and optimization iterations.

\section*{\bf{Graphical abstract}}
{\centering
\includegraphics[width=0.82\textwidth]{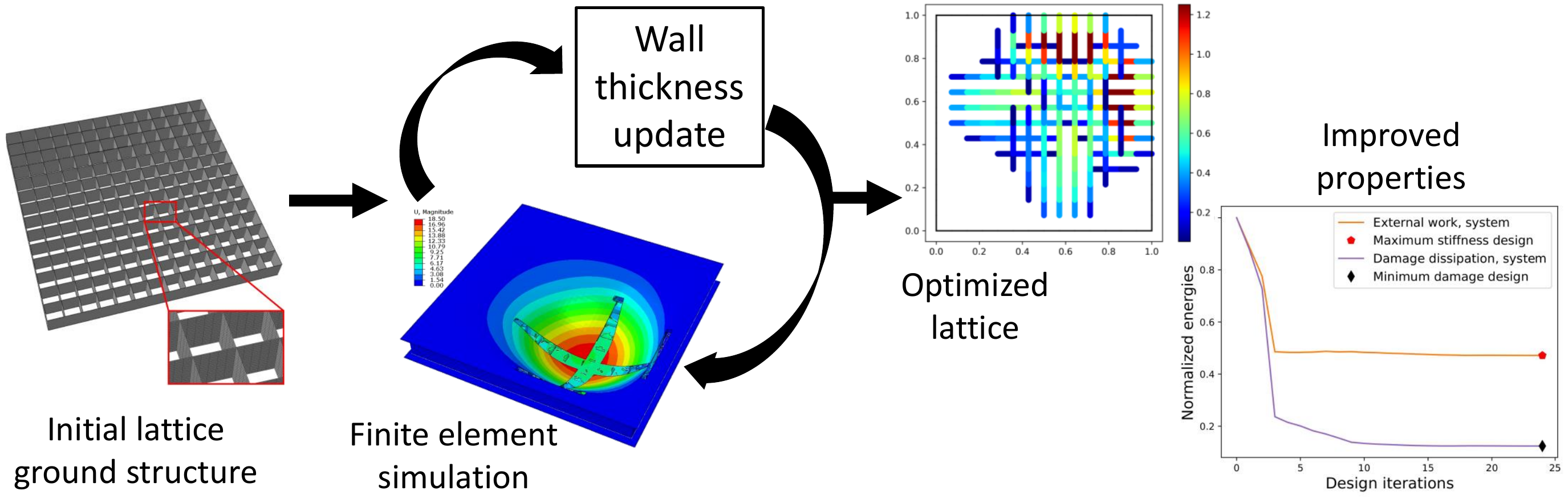}
\par
}
\end{abstract}

\begin{keyword}
Topology optimization \sep Thin-walled extruded lattice \sep Blast loading \sep High strain rate loading
\end{keyword}

\end{frontmatter}

\section{Introduction}
\label{sec:intro}
Thin-walled lattices are often employed to strengthen structures for energy absorbing purposes \citep{sun2010crashworthiness,yin2014crashworthiness,san2020thin}. Typically, these structures absorb energy by undergoing highly nonlinear large deformation that extends well beyond the linear elastic region \cite{duddeck2016topology}. Irreversible processes such as plastic deformation, damage, and fracture become the primary means of energy dissipation \cite{zhu2020energy,baroutaji2017crashworthiness}. It is of significant practical interest to leverage topology optimization (TO) algorithms to identify the optimum cross-sectional design and location of reinforcements \cite{duddeck2016topology}. The objectives are typically to maximize specific energy absorption \cite{yin2014crashworthiness,fu2019design,sun2014crashing}, maximize crashing force efficiency\cite{zarei2007crashworthiness,sun2014crashing}, and minimize peak crashing force \cite{wu2016design,abramowicz2003thin}. However, TO frameworks that are based on quasi-static loading, linear elastic, small-strain finite element (FE) formulations are unsuitable in this case, due to the large deformation and dynamic nature of the loading. For those that do account for geometric nonlinearility \cite{jung2004topology,clausen2015topology,wallin2016topology}, material plasticity \cite{abueidda2021topology,wallin2016topology,maute1998adaptive}, and damage \cite{verbart2016damage,li2017topology,james2015topology}, an in-house implementation of the FE procedure is required, as information on the material tangent stiffness matrix and values of the internal state variables are needed in the adjoint method to compute the gradient of the objective function \cite{tsay1990nonlinear,alberdi2018unified,abueidda2021topology}. This requirement is cumbersome for two reasons: (1) many general purpose FE codes, either open source or proprietary, already exist that can robustly solve large deformation nonlinear problems, they should be leveraged whenever possible; and (2) information on the tangent stiffness matrix may not be provided by those programs, or in the case of explicit dynamic simulations it is simply not formed. In particular, it is also noted that adjoint analysis for transient problems is complex and difficult to implement \cite{farrell2013automated} and can be as expensive as the forward FE analysis \cite{sigmund2011usefulness}.

Given that existing crashworthiness research widely employs commercial explicit dynamics FE code such as Abaqus/Explicit \cite{Abaqus2021} and LS-DYNA \cite{murray2007users}, it is highly relevant and of great importance to study how these codes can be coupled to TO frameworks that do not require information beyond what the commercial code provides. In the case of TO for thin-walled structures, several frameworks have been proposed and demonstrated. Notable ones are: (1) the equivalent static load approach, which seeks a static load that generates equivalent deformation to the true dynamic load, and thus converts the problem to a quasi-static TO problem \cite{lee2015nonlinear,jang2012dynamic}; (2) cellular automata-based approaches, which use a series of update rules to iteratively update the thickness of the shell elements defining the lattice walls \cite{hunkeler2013topology,duddeck2016topology,hunkeler2014topology,zeng2017improved}; (3) the response surface method, which seeks to approximate the true response surface of the objective function through simple functions and repeated FE simulations \cite{kurtaran2002crashworthiness,avalle2002design}; and (4) various probabilistic and evolutionary methods like Bayesian optimization \cite{liu2019design}, ant colony method \cite{liu2021multisurrogate,liu2019topographical}, and particle swarm method \cite{gao2019application}. Most of these methods, although capable of improving the design, take many FE simulation and optimization iterations to do so \cite{hunkeler2013topology}, and can be computationally expensive for large scale problems. Therefore, the current work aims at proposing a novel heuristic framework named LatticeOPT that can work with commercial FE code to generate lattice cross-section designs with improved specific energy absorption or minimized damage, with as few design iterations as possible. 

This paper is organized as follows: \sref{sec:methods} presents an overview of the definition of the lattice design space, thickness update schemes, and the workflow of the LatticeOPT algorithm. \sref{sec:results} presents and discusses the results from three numerical examples. \sref{sec:conc} summarizes the outcomes and highlights possible future works.

\section{Methods}
\label{sec:methods}
\subsection{Defining the lattice design space and design variables}
\label{Design_space}
\begin{figure}[h!] 
    \centering
     \subfloat[]{
         \includegraphics[width=0.55\textwidth]{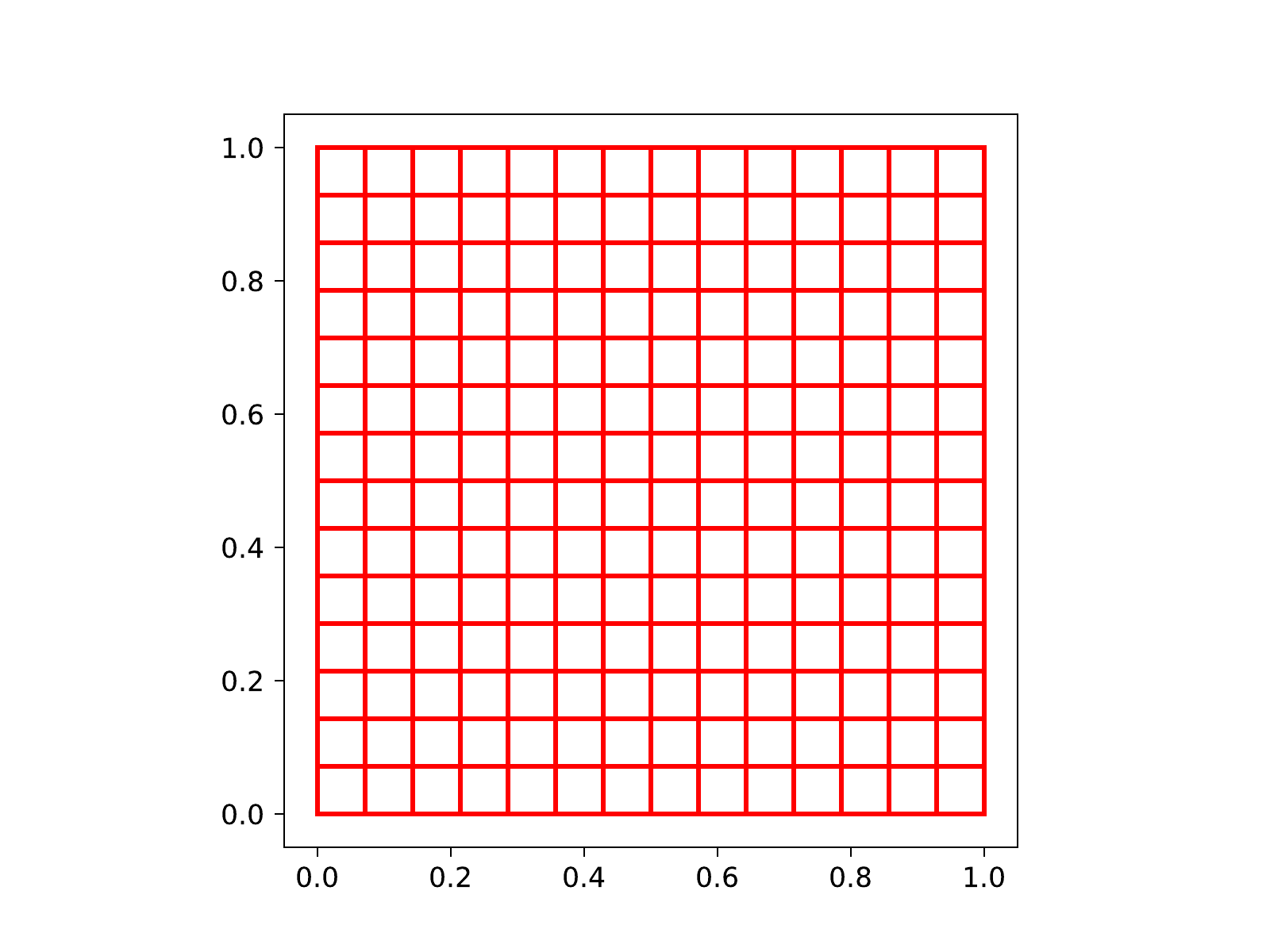}
         \label{fig:ground_structure}
     }
     \subfloat[]{
         \includegraphics[width=0.36\textwidth]{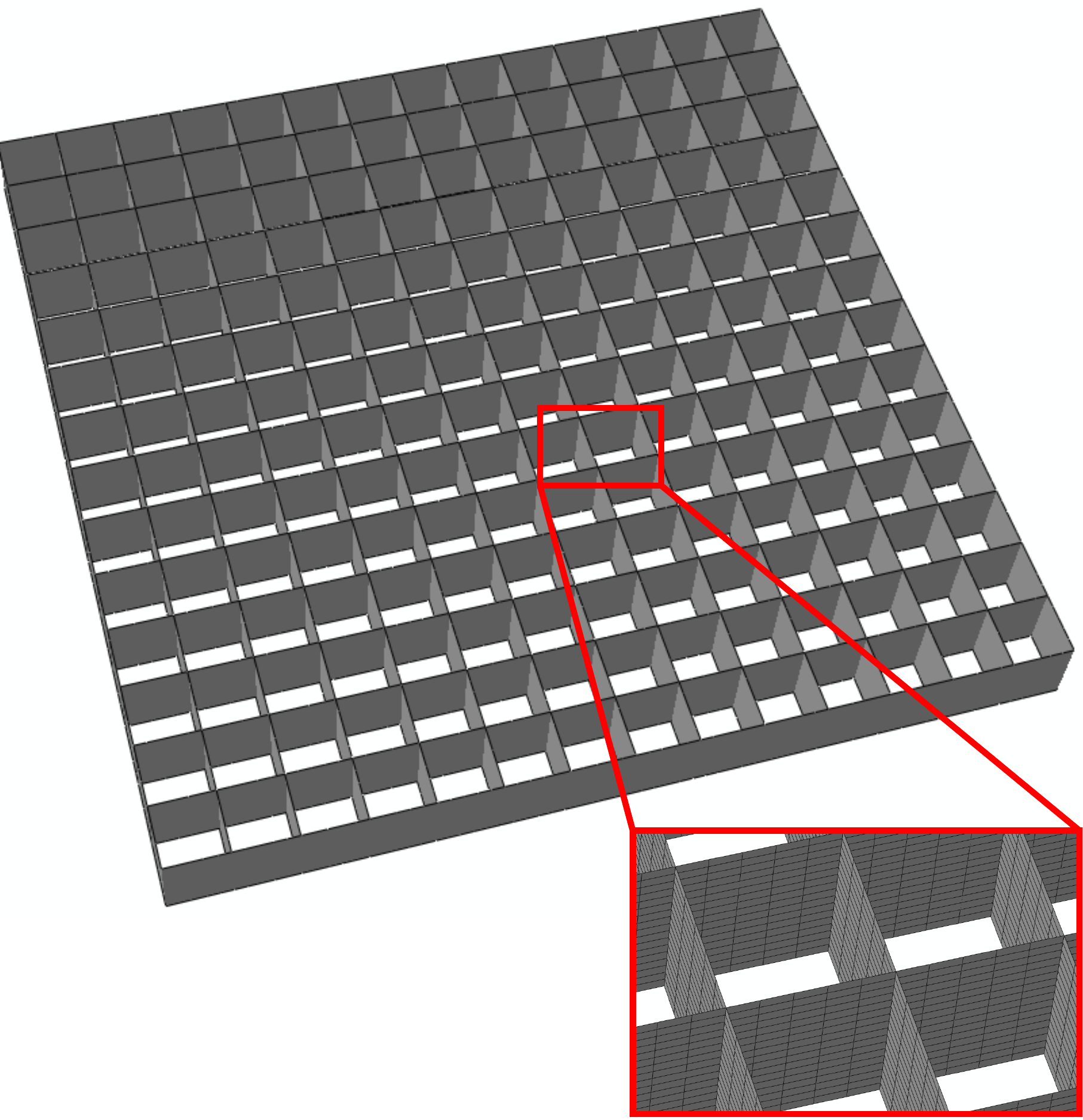}
         \label{fig:mesh}
     }
    \caption{Design space definition: \psubref{fig:ground_structure} The cross section is partitioned into a user-defined number of rectangular cells. Cells do not need to be squares. \psubref{fig:mesh} An extruded lattice is constructed from the cross section and a user-defined height, the inset shows the finite element mesh. Note that more than one element is used to discretize a single cell wall.}
    \label{initial_mesh}
\end{figure}

Currently, the LatticeOPT framework supports the definition of a cubic lattice design space, defined by the in-plane cross-section length (L) and width (W), as well as the out-of-plane height (H). Inspired by the ground structure approach \cite{soto1999basic}, which is typically applied to TO of truss structures, we partition the in-plane cross section into smaller rectangular cells, where each cell wall is represented by shell elements in FE simulations, see \fref{initial_mesh}. The user provides the number of cells in the X- and Y-directions, denoted by $N_x$ and $N_y$, respectively. A larger $N$ means smaller partitions and larger design space with more degrees of freedom. While in the Z-direction, the user can specify the number of finite elements used to discretize the lattice height, denoted by $N_z$.

Contrary to TO applied to solid structures, where the design variables are material densities of the elements, the design variables in LatticeOPT are the thicknesses of lattice walls, and the material property is not interpolated as is typically done in typical TO. In addition, typical TO frameworks penalize material densities using either the SIMP or RAMP schemes to avoid the presence of intermediate density values \cite{rozvany2009critical}, while in the present framework, no penalty scheme is needed. Wall thickness can be regarded as a continuous variable, and any intermediate thickness values within a reasonable range $[ t_{min} , t_{max} ]$ is physically meaningful and is manufacturable via additive manufacturing techniques such as selective laser melting \cite{frazier2014metal}. We would like to emphasize that, although the process of varying the shell element thickness is similar to a sizing optimization and/or topometry optimization \cite{leiva2004topometry}, since LatticeOPT usually involves many more design variables and it alters the cross section topology directly by removing edges below a minimum thickness, this is still considered to be a topology optimization framework. The use of the term topology optimization is also consistent with many previous studies on optimization of thin-walled shell structure cross sections that focus on updating shell thickness and remove cell walls below a thickness threshold \cite{hunkeler2013topology,hunkeler2014topology,duddeck2016topology}

\subsection{Thickness update schemes}
\label{update_schemes}
Many previous heuristic TO frameworks for thin-walled structures, like the hybrid cellular automata (HCA)-based framework by Hunkeler \cite{hunkeler2014topology} and Duddeck et al. \cite{duddeck2016topology}, evolve the design by iteratively updating the lattice wall thicknesses through some predetermined, intuition-based updating rules. The central idea embedded in these rules is the equalization of energy density: for efficient energy absorption, different parts of the lattice should contribute equally to the energy absorption \cite{hunkeler2014topology}. In terms of the ground lattice structure depicted in \fref{fig:mesh}, it means that each lattice wall should contribute equally to the total energy absorption \cite{hunkeler2014topology}. 

Since a design can be uniquely characterized by the wall thickness of each lattice wall, design update boils down to the update of lattice wall thicknesses. Inspired by the idea of equalization of energy density, we propose two different thickness update schemes, which produce the updated thickness distribution $\bm{t}^{i+1}$ for the next design iteration.

\subsubsection{Scheme 1}
\label{std_scheme}
Given the current wall thickness distribution $\bm{t}^i$ (superscript $i$ denotes the design iteration index), FE simulation provides the energy $\bm{E}^i$ absorbed by each lattice wall, which accounts for the stored elastic strain energy and energy dissipated due to plasticity and damage. Approaching equalization of energy density at the next design iteration is mathematically equivalent to the following optimization problem:

\begin{equation}
    \min_{ \forall {\bf{A}}^{i+1} \ge 0 } \; {\rm{std}}( \frac{ {\bf{E}}^{i+1} }{ {\bf{A}}^{i+1} } ) \footnote{Note that we have replaced the lattice wall volume with its projection area, since all walls share a constant height.} \, ,
    \label{min_std_1}
\end{equation}
where the projection area of lattice wall $j$ at design iteration $i+1$ (denoted by $A^{i+1}_j$) is given by $A^{i+1}_j = L_j t^{i+1}_j$, with $L_j$ being the length of that lattice wall. In the case of a uniform energy density distribution, the standard deviation attains its global minimum of 0. However, \eref{min_std_1} offers little value in practice: $\bm{E}^{i+1}$ is not known \textit{a priori}, nor can it be estimated to first order using a truncated Taylor series given a thickness change $\Delta \bm{t}$:
\begin{equation}
    {\bf{E}}^{i+1} \approx {\bf{E}}^{i} + \frac{ d \bm{E}^{i} }{ d \bm{t}^{i} } \Delta \bm{t} \, .
    \label{Taylor_estimate}
\end{equation}
In general, it is difficult to obtain the gradient $\frac{ d \bm{E}^{i} }{ d \bm{t}^{i} }$ through the adjoint method since the FE tangent stiffness matrix (which is needed in the calculation of the adjoint vector) is not always available from commercial FE packages, and for explicit dynamic simulations, the tangent stiffness matrix is simply not formed. The sheer size of the gradient matrix also makes it prohibitively expensive to estimate using the finite difference method. Therefore, we modify the optimization problem using a zeroth order estimate for $\bm{E}^{i+1}$:

\begin{equation}
    \min_{ \forall {\bf{A}}^{i+1} \ge 0 } \; {\rm{std}}( \frac{ {\bf{E}}^{i} }{ {\bf{A}}^{i+1} } ) \, .
    \label{scheme1}
\end{equation}
Note that now the only unknowns are the updated thickness array $\bm{t}^{i+1}$ embedded inside the projection area array $\bm{A}^{i+1}$, and the gradient of the objective function is straightforward ($\bm{E}^{i}$ is independent of $\bm{t}^{i+1}$). This is beneficial, as gradient-based optimizers can be employed to efficiently reduce the objective function, which is highly urged in the famous forum paper \cite{sigmund2011usefulness}. For practical purposes, we further require the following variable bounds, volume, and thickness change constraints:
\begin{equation}
    t^{i+1}_j \in [ 0 , t_{max} ], \forall j = 0 , 1 , ... , n
    \label{t_bnd}
\end{equation}
\begin{equation}
    H \sum{ {\bf{A}}^{i+1} } - V^* = 0 ,
    \label{mass_con}
\end{equation}
\begin{equation}
    L_{\infty} ( {\bf{t}}^{i+1}  - {\bf{t}}^{i}  ) \le \Delta t_{max} ,
    \label{dt_con}
\end{equation}
where $V^*$ and $H$ denote the target volume and lattice height, respectively. In practice, the maximum thickness change constraint (\eref{dt_con}) is converted to bounds on the design variables $\bm{t}$ during the optimization. \eref{scheme1} - \eref{dt_con} define a nonlinear constrained optimization problem, the solution to which provides the updated thickness distribution $\bm{t}^{i+1}$.

\subsubsection{Scheme 2}
\label{beso_scheme}
The second updating scheme is inspired by the bi-directional evolutionary structural optimization (BESO) method, pioneered by Huang and Xie \cite{huang2007convergent}. It is a heuristic TO method typically applied to solid structures, where each element in the mesh is ranked by a so-called sensitivity number, and elements are added or removed based on their relative sensitivity ranking. When optimizing for SEA, the element sensitivity $\alpha$ at iteration $i$ is defined as \cite{doi:10.1080/13588260701497862}:
\begin{equation}
    \bm{\alpha}^{i} = \frac{ \bm{V} }{ V_{tot} } - \frac{ \bm{E}^{i} }{ E_{tot}^{i} },
    \label{beso_sen}
\end{equation}
where $\bm{V}$ and $\bm{E}$ denote the undeformed volume and absorbed energy of the elements, respectively. Following this definition, we can analogously define the sensitivity for lattice walls at design iteration $i+1$ as:
\begin{equation}
    \bm{\alpha}^{i+1} = \frac{ \bm{A}^{i+1} }{ \sum{ {\bf{A}}^{i+1} } } - \frac{ \bm{E}^{i+1} }{ E_{tot}^{i+1} } \, .
    \label{beso_sen_wall}
\end{equation}

Equalization of energy density can again be expressed as minimizing the standard deviation of the lattice wall sensitivity number. Estimating $\bm{E}^{i+1}$ with $\bm{E}^{i}$, the second updating scheme can be stated in terms of the following optimization problem:

\begin{equation}
    \min_{ \forall {\bf{A}}^{i+1} \ge 0 } \; {\rm{std}}( \frac{ \bm{A}^{i+1} }{ \sum{ {\bf{A}}^{i+1} } } - \frac{ \bm{E}^{i} }{ E_{tot}^{i} } ) \, ,
    \label{scheme2}
\end{equation}
subject to the same bounds and constraints defined in \eref{t_bnd} - \eref{dt_con}.

\subsubsection{Implication of the update schemes}
\label{implication}
Note that both update schemes introduced above originate from optimizing energy absorption of the lattice. However, it is not immediately obvious how the attempted equalization of energy density actually affects the performance of the updated lattice.

As we shall see in \sref{sec:results}, the two schemes actually lead to improved stiffness of the optimized lattice structure. In the case of a fixed-displacement loading scenario (e.g., \sref{sec:ex1}), a stiffer structure absorbs more energy compared to its more compliant counterparts, thus agreeing with the expected objective of maximizing SEA. We also demonstrate in the case of a fixed-load loading scenario (e.g., \sref{sec:ex2}) that, instead of maximizing energy absorption, a stiffer structure minimizes material damage.

\subsection{Periodic lattice structures}
\label{Periodic}
LatticeOPT supports the generation of periodic lattices. The user can enable this functionality by specifying the number of unit cells in the X- and Y-directions, denoted by $N^{uc}_x$ and $N^{uc}_y$. The algorithmic treatment is inspired by the work of Huang and Xie \cite{huang2008optimal}, where the BESO method is applied to generate periodic lattice structures. In this case, the effective design variables become the thickness of all lattice walls in the periodic unit cell, denoted as $\Tilde{\bm{t}}$. The absorbed energy of each lattice wall in the system is aggregated to the periodic unit cell through the following aggregation matrix $\bm{A}$:

\begin{equation}
\underset{ n_{uc} \times 1 }{  \begin{bmatrix}
    \Tilde{E}_1 \\
    \Tilde{E}_2 \\
    \vdots \\
    \Tilde{E}_{n_{uc}} \\
    \end{bmatrix} } =  \underset{ \bm{A} \, , \, n_{uc} \times n_{f} }{ \begin{bmatrix}
    1 & 0 & 0 & 0 &  \cdot \cdot \cdot  & 1 & 0 & 0 & 0  \\
    0 & 1 & 0 & 0 &  \cdot \cdot \cdot  & 0 & 1 & 0 & 0  \\
    \vdots & & & & \ddots & & & & \vdots \\
    0 & 0 & 0 & 1 & \cdot \cdot \cdot  & 0 & 0 & 0 & 1
    \end{bmatrix} }  \cdot
\underset{ n_{f} \times 1 }{ \begin{bmatrix}
    E_1 \\
    E_2 \\
    \vdots \\
    \vdots \\
    E_{n_f} \\
    \end{bmatrix} },
\end{equation}
where $n_{uc}$ and $n_f$ denote the number of lattice walls in the unit cell and full design space. The nonzero entries in the aggregation matrix indicate which wall entry in the unit cell does the current wall correspond to, and the aggregation operation effectively sums the absorbed energies of all periodic images of a particular wall in the unit cell. The aggregated energy $\bm{\Tilde{E}}^{i}$ is fed into the thickness update scheme to generate $\Tilde{\bm{t}}^{i+1}$. The complete thickness array $\bm{t}^{i+1}$ can be recovered by the scattering operation:
\begin{equation}
    \underset{ n_{f} \times 1 }{ \bm{t}^{i+1} } = \underset{ n_{f} \times n_{uc} }{ \bm{A}^T } \cdot
    \underset{ n_{uc} \times 1 }{ \Tilde{\bm{t}}^{i+1} },
\end{equation}
which enforces that all periodic images of a particular wall in the unit cell must share identical thickness.

\subsection{Implementation}
\label{implementation}
The workflow of the LatticeOPT framework is summarized in Alg. (\ref{LatticeOPT}). To begin the analysis, the user supplies the design space dimensions and the number of periodic unit cells. The user also needs to provide information about the material properties, loading, and BCs in the form of a FE input file compatible with the analysis package of choice. For this work, we have selected Abaqus \cite{Abaqus2021} to be the analysis tool, but we emphasize that the LatticeOPT framework can work with any other FE packages so long as the user provides a compatible FE input file. A linear thickness scaling is applied to the initial thickness distribution to satisfy the volume constraint directly (lines 1-2). The aggregation matrix is formed if the system contains more than one unit cell (lines 3-4). A uniform mesh of the ground lattice structure containing all lattice walls in the system, is generated and added to the FE input file (line 5). 

With all the preparation work done, we now begin the optimization. The current thickness distribution is added to the input file, from which a FE simulation is conducted, and the wall energy array is calculated (lines 8-10). Two different performance metrics are calculated for the current design, they are (1) SEA:
\begin{equation}
    SEA = \frac{ \sum{ \bm{E}^i } }{ V^* },
\end{equation}
and in the presence of material damage (2) mean wall connectedness (MWC):
\begin{equation}
    MWC = \underset{ 1 \le j \le N_{wall} }{ {\rm{mean}} } \frac{ N_{connect,j} }{ N_z },
    \label{eq:MWC}
\end{equation}
where $N_{connect,j}$ denotes the number of layers of finite elements (out of a total of $N_z$ layers) in the Z-direction that remain \emph{totally undamaged} for lattice wall $j$ in the design space. MWC provides a quantitative measure of how much the lattice walls remain interconnected after the applied loading. MWC falls in the range $[0,1]$. A lower MWC means that more lattice walls are damaged and become disconnected from their neighbours. Both SEA and MWC are evaluated at the end of the load step.

Following the recommendation in \cite{huang2008optimal}, we average the wall energy array in the current iteration with that from the previous iteration to stabilize the evolution (line 12). Note that unlike classical TO, a mesh independence filter is not applied to the energy array, since each lattice wall contains multiple finite elements, and $E^i_j$ sums the energy contribution of all those elements. Thus $\bm{E}$ is expected to be less prone to mesh dependence. If the structure is periodic, wall energies are aggregated to the unit cell (line 14). Before solving the optimization problem, we set the variables' bounds based on their current values and the maximum allowable thickness change per iteration, as well as the range of allowable thicknesses (line 16). To solve the nonlinear optimization problem resulted from the thickness updating scheme (line 17), we used the \emph{minimize} function from Scipy \cite{2020SciPy-NMeth}, and in particular, the trust region constrained algorithm (trust-constr) \cite{conn2000trust,byrd1987trust} is used. It is worth noting that different nonlinear optimizers were tested, such as the sequential least-squares programming (SLSQP), the interior point method (IPOPT) \cite{biegler2009large} and the method of moving asymptotes (MMA) \cite{svanberg1987method}. It was found that IPOPT and trust-constr gave similar computational performance and were more robust than other methods tested. The trust-constr method was selected for its performance and the advantage that it is a built-in package in Scipy, which is a widely used Python package. Finally, periodic scattering is applied to recover the full thickness distribution if the system is periodic (line 19). 

To ensure that the generated designs are manufacturable, we remove walls whose thickness fall below the threshold $t_{min}$ (line 20). Similar to the soft-kill BESO method \cite{yang1999bidirectional}, we assign a small thickness ($10^{-6}$ mm) to those 'removed' walls, instead of actually removing them from the FE model. After wall removal, a linear thickness scaling is again applied to ensure that the thresholded design also satisfies the volume constraint (line 22). Finally, the algorithm checks for convergence by monitoring the maximum thickness change (line 24), and moves on to the next design iteration if convergence is not achieved (line 26).

\begin{algorithm}[!ht]
\DontPrintSemicolon
    \KwInput{ $[N_x,N_y,N_z]$ , $[L,W,H]$ , $[N^{uc}_x,N^{uc}_y]$ , loading, BCs , material properties, max design iteration , Update scheme , $[ \, \bm{t}^0 , t_{min} , t_{max} \, ]$ , $V^*$ , $\Delta t _{max}$}
    
    \KwOutput{Optimized design defined by thickness distribution array $\bm{t}^*$ }

    \tcc{Initialization}
    $ r \gets \frac{ V^* }{ \sum{ \bm{A}^0} H }$ \tcp*{Initial volume ratio}
    $ \bm{t}^0 \gets \bm{t}^0 * r $ \tcp*{Scaling to satisfy volume constraint}
    \If{ $N^{uc}_x N^{uc}_y > 1$ }
    {
        Form aggregation matrix $\bm{A}$ \tcp*{Only need to form once}
    }
    Write FE input file for ground lattice structure
    
    $ i \gets 0 $
    
    \tcc{Begin optimization}
    \While{ $i < $ max design iteration }{
    Add $\bm{t}^i$ information to FE input file
    
    Run FE simulation, obtain $\bm{E}^i$
    
    Compute performance metrics
    
  \If{ $ i > 0 $ }
    {
        $ \bm{E}^i \gets \frac{1}{2} ( \bm{E}^{i-1} + \bm{E}^i )$ \tcp*{Average energy with last iteration}
    }

    \If{ $N^{uc}_x N^{uc}_y > 1$ }
    {
        $\Tilde{\bm{E}^i} = \bm{A}  \bm{E}^i $ \tcp*{Periodic aggregate}
    }
    
    \For{ $ j = [ 1 : N_{walls} ] $ }{
    \tcc{Set variable bounds}
    $bnd_j = [ \, max( 0 , t^i_j - \Delta t _{max} ) \, , \, min( t^i_j + \Delta t _{max} , t _{max} ) \, ]$
    }
    
    Solve optimization problem using optimizer to update thickness

    \If{ $N^{uc}_x N^{uc}_y > 1$ }
    {
        $\bm{t}^{i+1} = \bm{A}^T  \Tilde{\bm{t}}^{i+1}  $ \tcp*{Periodic scatter}
    }

    $\bm{t}^{i+1} [ \;  \bm{t}^{i+1} < t_{min} \; ] = 10^{-6}$ \tcp*{Thresholding}
    
    $ r \gets \frac{ V^* }{ \sum{ \bm{A}^{i+1}} H }$ \tcp*{Current volume ratio}
    
    $ \bm{t}^{i+1} \gets \bm{t}^{i+1} * r $ \tcp*{Scaling to satisfy volume constraint}
    
    Save $\bm{E}^i$, $\bm{t}^{i+1}$ and performance metrics to file
    
    \tcc{Check convergence}
  \If{ $ max( \bm{t}^{i+1} - \bm{t}^{i} ) < 10^{-2} $ }
    {
        break
    }
    
    $ i \gets i + 1 $ \tcp*{Move to next design iteration}

    }
\caption{LatticeOPT}
\label{LatticeOPT}
\end{algorithm}

\section{Results and discussion}
\label{sec:results}
In this section, we present three numerical examples that showcase the capabilities of LatticeOPT. In the first example, we performed a benchmark test against the HCA framework developed by Hunkeler \cite{hunkeler2014topology}, as well as the commercial TO package LS-OPT \cite{goel2009topology} studied therein. In the second example, we performed optimization of a square lattice sandwiched between two face sheets, subjected to a blast loading, to maximize its stiffness and minimize damage. In the last example, we repeated example 2, but enforced a 2-by-2 periodic unit cell arrangement. All simulations presented in this section were conducted in Abaqus/Explicit \cite{Abaqus2021} on an Intel i7-11800H processor using 8 cores.

\subsection{Benchmark test against HCA and LS-OPT}
\label{sec:ex1}
Hunkeler \cite{hunkeler2014topology} studied the TO of a long, slender column under dynamic axial compression using a HCA-based algorithm. The column has a cross sectional dimensions of 80-by-100 mm$^2$, and a height of 400 mm. The out-most rectangular boundaries of the column are nondesignable and have a fixed wall thickness of 1.5 mm. The cross section is partitioned into 4 cells in the X direction and 5 cells in the Y direction, leading to a total of 31 designable lattice walls. However, symmetry boundary conditions were used by Hunkeler \cite{hunkeler2014topology}, so only a quarter of the cross section with a total of 10 designable edges were included. The material was aluminum with elasto-plastic material behavior modeled with a rate-independent Mises plasticity model with a piece-wise linear isotropic hardening curve. The hardening curve is provided in \tref{flow_stress}, which was taken from the data in \cite{hunkeler2014topology} (see Table 5.4 therein). The column was fixed at the bottom, and was subjected to compression by a rigid plate having an initial downward velocity of 5 m/s and a mass of 500 kg. To ensure comparability of results, mesh sizes identical to that in \cite{hunkeler2014topology} were used to discretize the structure. 

\begin{table}[h!]
    \caption{Piece-wise linear hardening curve used in the work of Hunkeler \cite{hunkeler2014topology}}
    \small
    \centering
    \begin{tabular}{cccccccccc}
    \hline
    Flow stress (MPa) & \vline & 180.0  & 190.0 & 197.0 & 211.5 & 225.8 & 233.6 & 238.5 & 248.5 \\
    \hline
    Equivalent plastic strain & \vline & 0.0  & 0.01 & 0.02 & 0.05 & 0.1 & 0.15 & 0.2 & 0.4 \\
    \hline
    \end{tabular}
    \label{flow_stress}
\end{table}

The original optimization problem studied in \cite{hunkeler2014topology} was to minimize mass of the structure while maintaining an end displacement of less than 75 mm after 25 ms. However, the author did report structure SEA at 10 ms for the optimized designs generated by HCA and LS-OPT (see Table 5.5 therein). As LatticeOPT focuses more on optimizing for energy absorption and stiffness enhancement instead of controlling end displacement, we reformulated the optimization problem to maximize SEA at 10 ms, subject to a structure mass constraint of 870 g. This target mass is comparable to that of the optimized designs from HCA and LS-OPT, which have masses of 878 and 863 g, respectively (see Table 5.5 therein). Two different TO runs were completed, one with each update schemes presented in \sref{update_schemes}, and the maximum allowable iterations was set to 25. The minimum and maximum allowable thicknesses were set to 0.4 mm and 2 mm, respectively.

To confirm comparability of results, the trivial design with only the out-most boundary walls present was simulated, and SEA at 10 ms was found to be 2.09 kJ/kg, very comparable to the 2.01 kJ/kg originally reported in \cite{hunkeler2014topology}. The small difference can be attributed to the use of symmetry boundary conditions in \cite{hunkeler2014topology}, which resulted in a slightly different deformation. 

\begin{figure}[h!] 
    \newcommand\x{0.85}
    \centering
         \includegraphics[trim={0 1cm 0cm 2cm},clip,width=\x\textwidth]{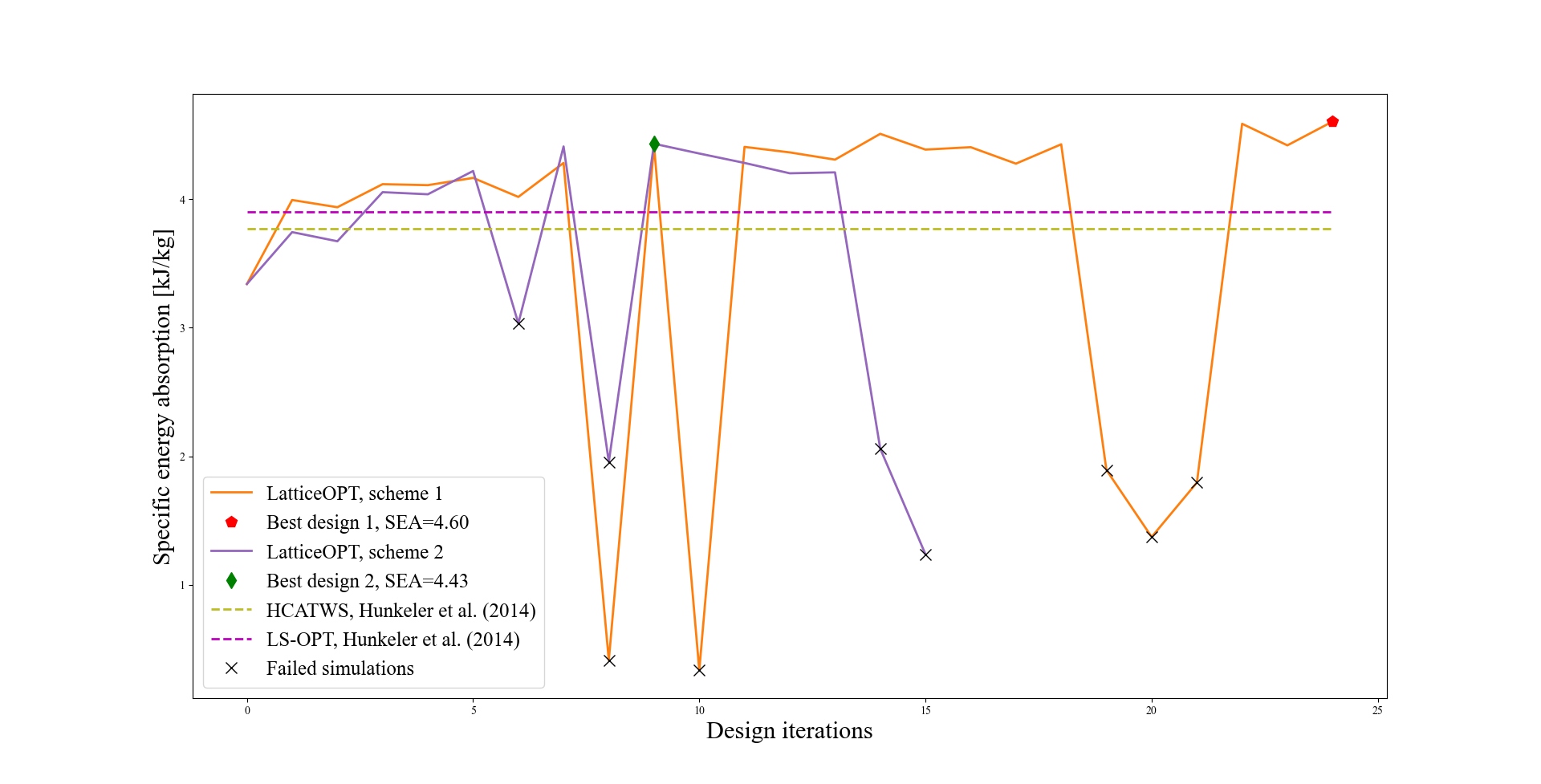}
    \caption{Evolution history of specific energy absorption. Note that both updating schemes yielded an optimized structure that has larger SEA than both the HCA and LS-OPT method.}
    \label{case1_comp}
\end{figure}

\begin{figure}[h!]
    \newcommand\x{0.29}
    \centering
    \begin{tabular}{ c  c  c  }
    \begin{minipage}[c]{\x\textwidth}
       \centering 
        \subfloat[Scheme 1, iteration 1]{\includegraphics[trim={2cm 0 2cm 0},clip,width=\textwidth]{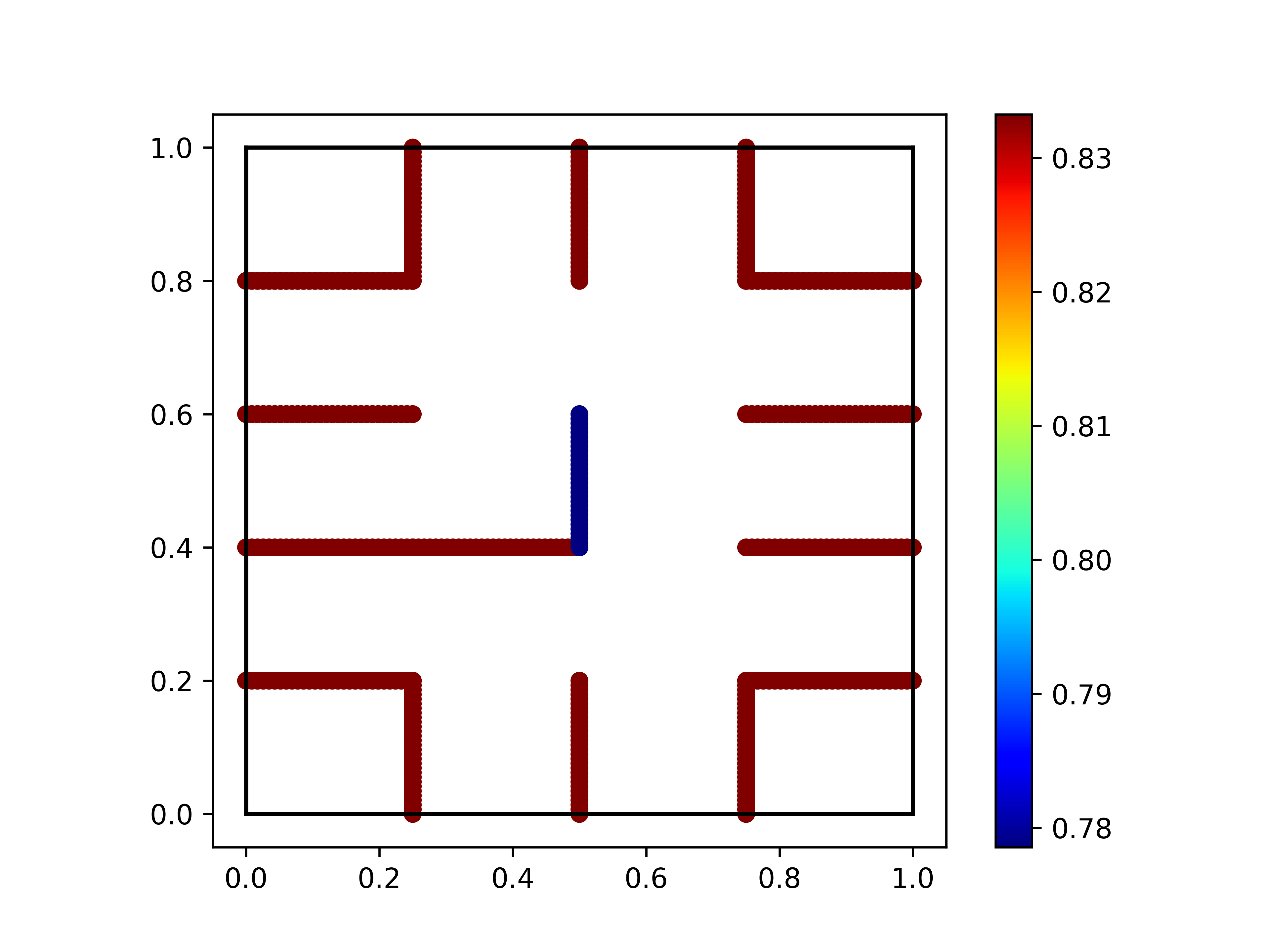}
        \label{fig:d1}}
    \end{minipage}
    &
    \begin{minipage}[c]{\x\textwidth}
       \centering 
        \subfloat[Scheme 1, iteration 12]{\includegraphics[trim={2cm 0 2cm 0},clip,width=\textwidth]{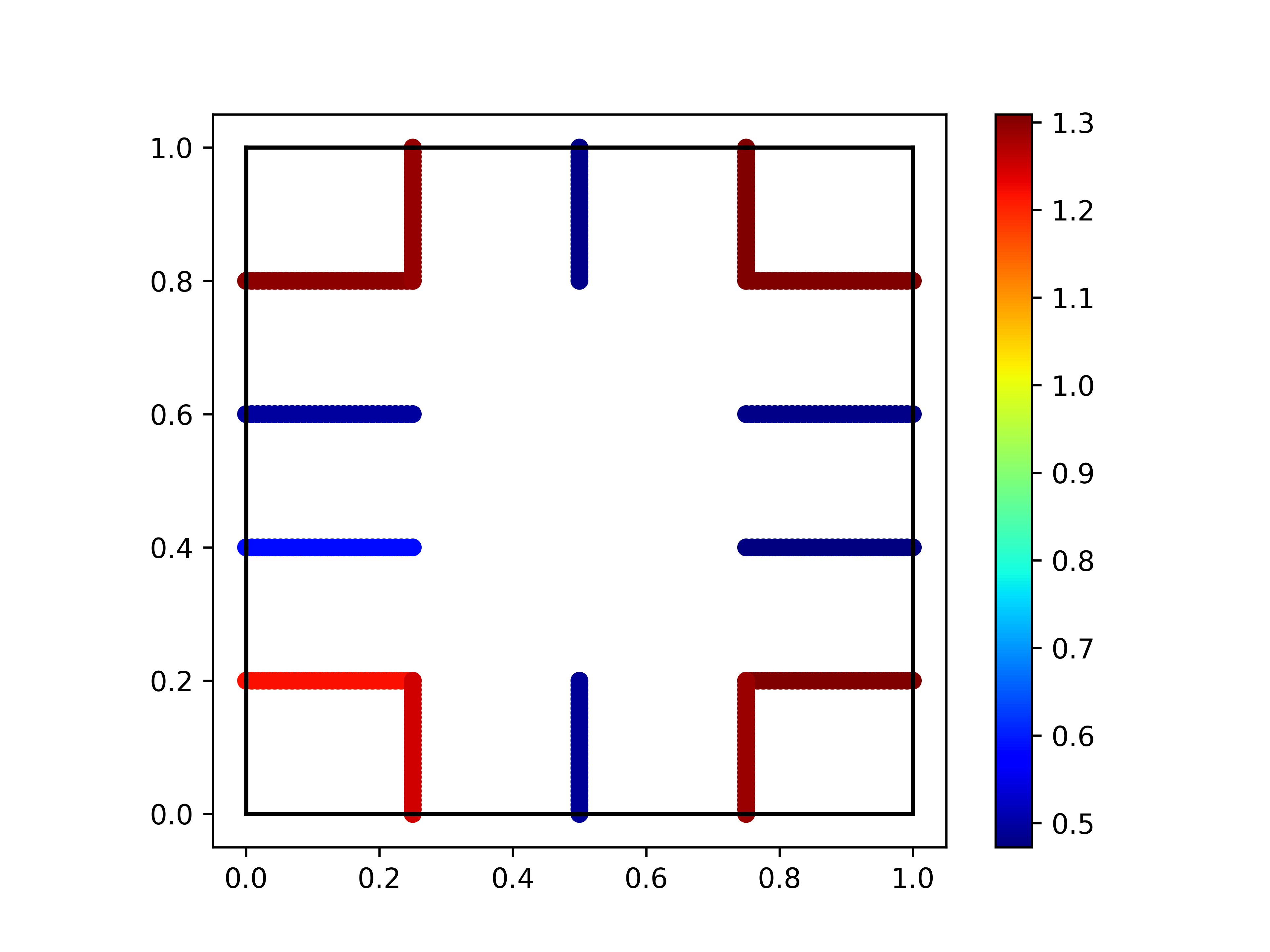}
        \label{fig:d2}}
    \end{minipage}
    &
    \begin{minipage}[c]{\x\textwidth}
       \centering 
        \subfloat[Scheme 1, iteration 24 $^*$]{\includegraphics[trim={2cm 0 2cm 0},clip,width=\textwidth]{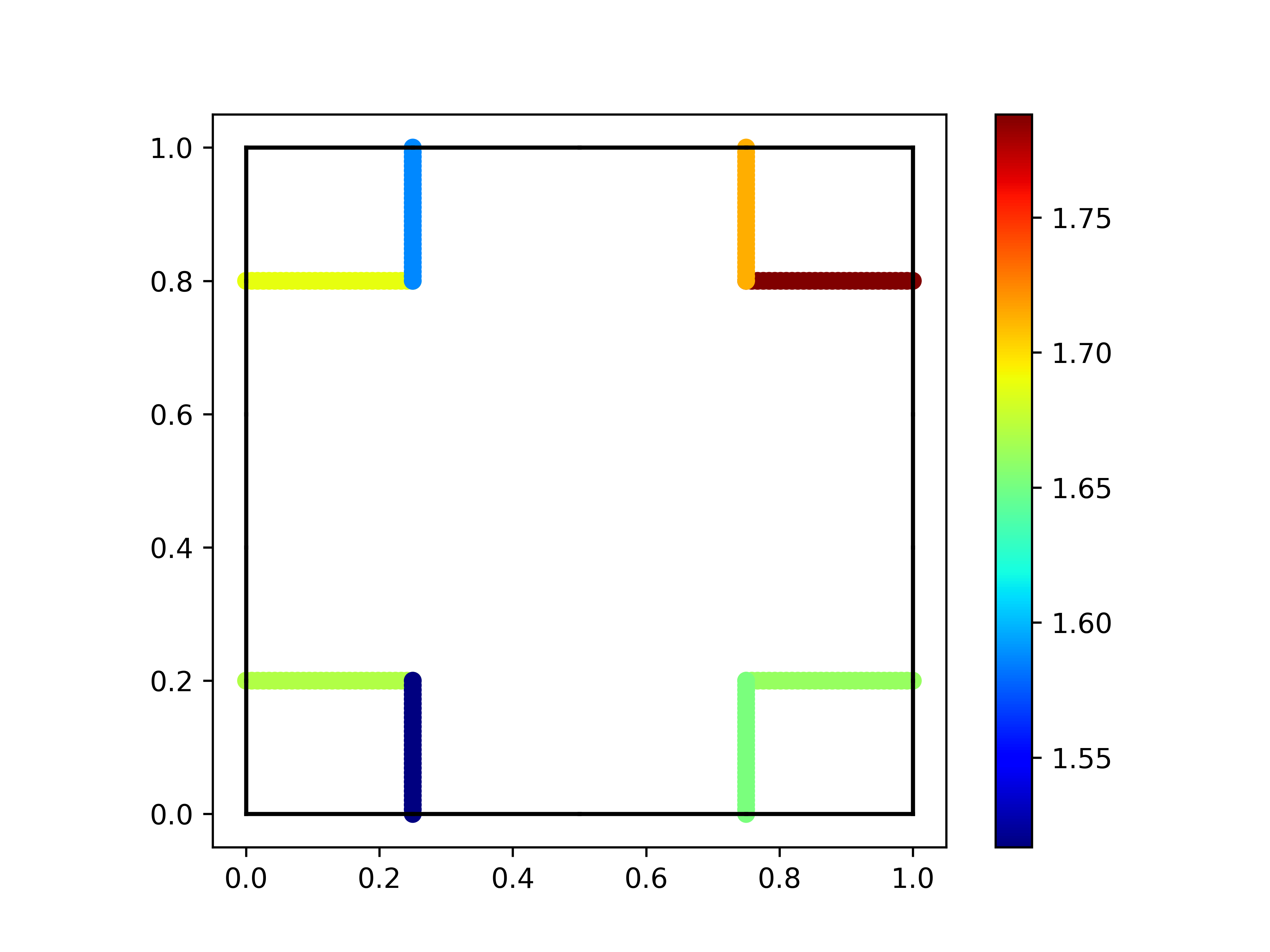}
        \label{fig:d3}}
    \end{minipage}\\

    \begin{minipage}[c]{\x\textwidth}
       \centering 
        \subfloat[Scheme 2, iteration 1]{\includegraphics[trim={2cm 0 2cm 0},clip,width=\textwidth]{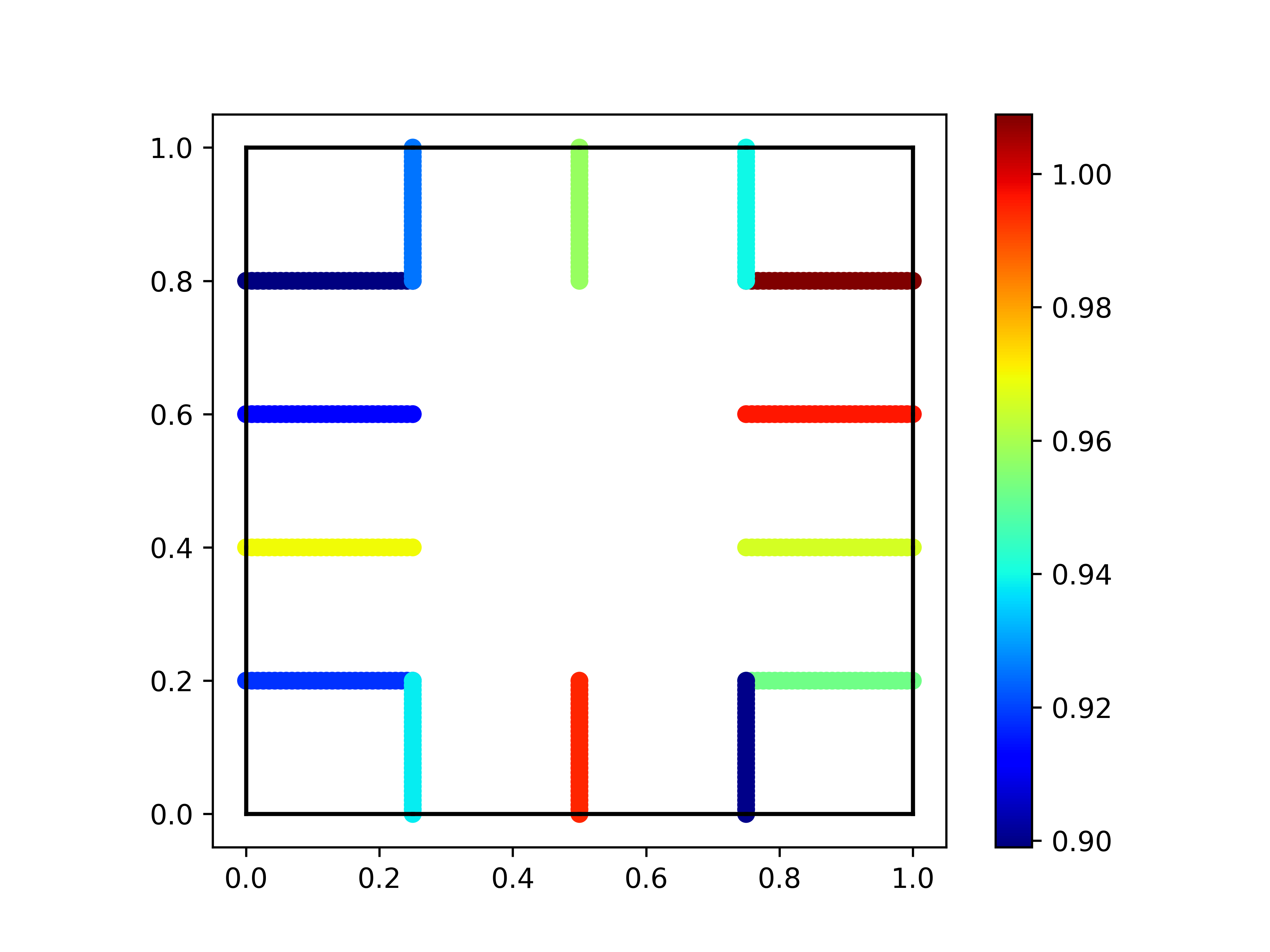}
        \label{fig:d4}}
    \end{minipage}
    &
    \begin{minipage}[c]{\x\textwidth}
       \centering 
        \subfloat[Scheme 2, iteration 9 $^*$]{\includegraphics[trim={2cm 0 2cm 0},clip,width=\textwidth]{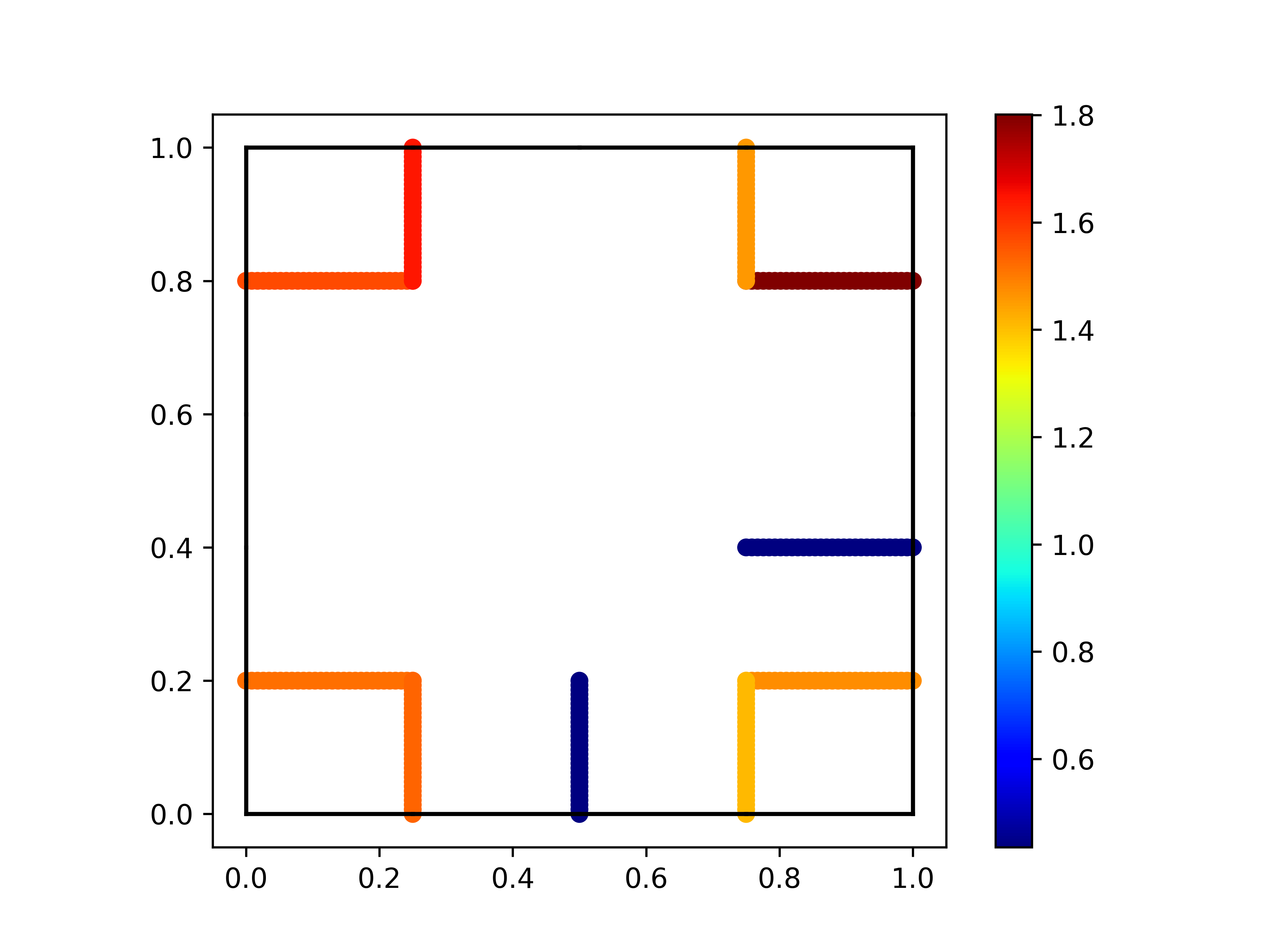}
        \label{fig:d5}}
    \end{minipage}
    &
    \begin{minipage}[c]{\x\textwidth}
       \centering 
        \subfloat[Scheme 2, iteration 15]{\includegraphics[trim={2cm 0 2cm 0},clip,width=\textwidth]{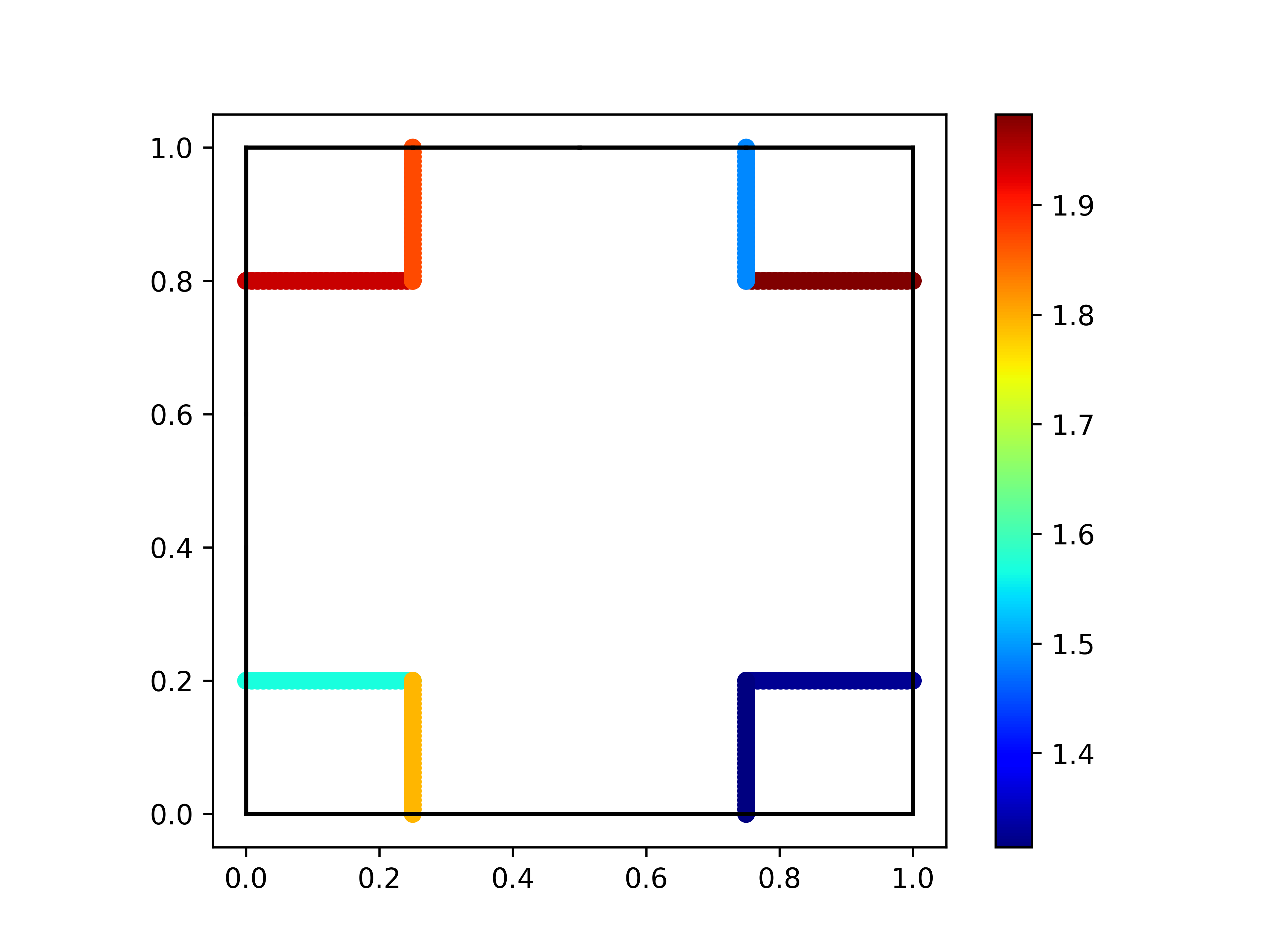}
        \label{fig:d6}}
    \end{minipage}
    \end{tabular}
    \caption{Designs generated by two update schemes at different design iterations, the optimized design that gives highest specific energy absorption is marked with the superscript $^*$. The color is assigned based on lattice wall thickness. The non-designable boundary walls are marked as black lines. }
    \label{case1_designs}
\end{figure}

When using a convergence tolerance of 0.01 mm (see line 24 of Alg. (\ref{LatticeOPT})), thickness update scheme 1 did not converge on a design and ran the full 25 design iterations, while update scheme 2 converged after 16 iterations. The SEA evolution histories for both update schemes, along with those from HCA and LS-OPT, are shown in \fref{case1_comp}. Note that the iteration indices in LatticeOPT are zero-based, so the initial design is iteration 0, and the last design iteration is 24 (for a total of 25 iterations). From the evolution history, we see that neither of the thickness update schemes guarantees increment in the objective function as more design iterations are conducted. This is reasonable, as the LatticeOPT framework is heuristic-based instead of gradient-based, so there is no strict guarantee that the objective function will improve every iteration. We also note that due to that large deformation nature of the problem, not all the simulations ran to completion (failed simulations are marked with black crosses in \fref{case1_comp}). LatticeOPT does not require the successful completion of the simulation to proceed. In the case of a failed simulation, $\bm{E}$ is simply extracted from the last output step in the simulation output, hence the lower SEA values for the failed iterations. Despite the oscillatory evolution history, we see that both update schemes produced an optimized design that has higher SEA than the best designs from HCA and LS-OPT. In this example, since the mass of the rigid impactor is much larger than the column, the impactor velocity was approximately constant, this loading case can be approximated as a fixed-displacement case. In this scenario, it is obvious from \fref{case1_comp} that the apparent consequence of the thickness update schemes is to maximize structure SEA, which implies that the structure stiffness is also improved. Selected intermediate designs generated by both schemes at different iterations are shown in \fref{case1_designs}. 

From \fref{case1_designs}, we see that both update schemes quickly removed lattice walls near the center of the column. As the optimization progressed, lattice walls on the four corners of the rectangular column were thickened, and this trend continued to the end of the optimization, where they became the only remaining lattice walls in the designable space. This design feature is consistent with the best designs generated from HCA (see Figure 5.35 therein) as well as LS-OPT (see Figure 5.41 therein). The trend where the corner walls get strengthened is also consistent with the findings described in the work of Kim et al. \citep{kim2002new}, which found this to be an effective strategy for increasing SEA. We also note that our generated designs are asymmetric due to the lack of symmetry boundary conditions, but this is considered more realistic as a long, slender column tends to buckle in compression, thus breaking the plane of symmetry. The larger design space (31 walls in our case vs. 10 in \cite{hunkeler2014topology}) is likely the reason for the higher SEA of the optimized designs generated from LatticeOPT. 

Finally, it is worth highlighting the effectiveness of LatticeOPT: both update schemes were able to generate high SEA designs with 25 or less FE simulations. As a comparison, the best HCA design was found after 84 simulations, and the best LS-OPT design after 187 simulations. The ability of LatticeOPT to generate optimized designs effectively is vital when the design space is large and results in an expensive FE simulation.

\subsection{Optimization of a lattice-filled sandwich panel under blast loading}
\label{sec:ex2}
\begin{figure}[h!] 
    \centering
     \subfloat[]{
         \includegraphics[width=0.5\textwidth]{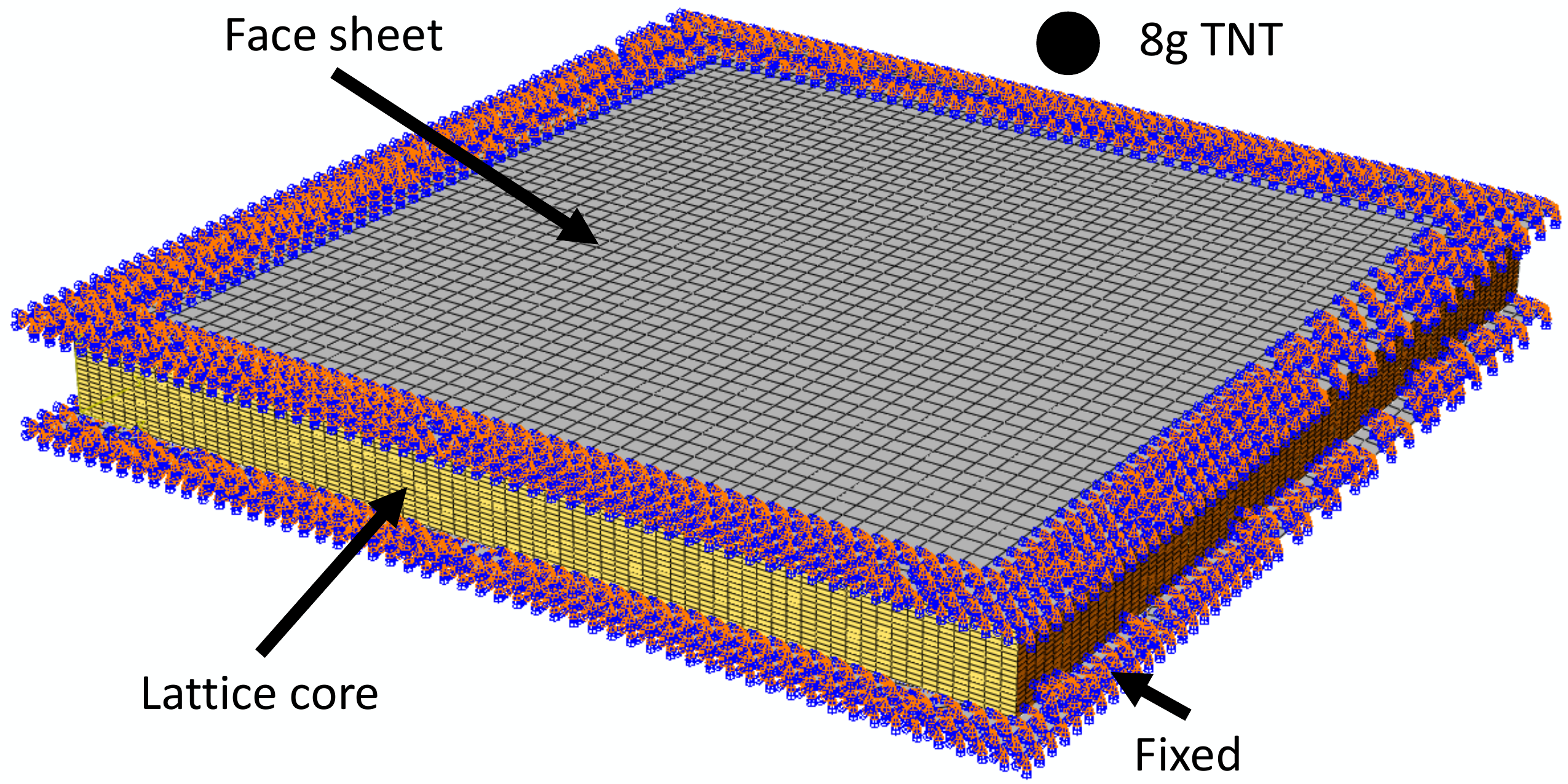}
         \label{fig:BCs}
     }
     \subfloat[]{
         \includegraphics[trim={0 0 18cm 0},clip,width=0.4\textwidth]{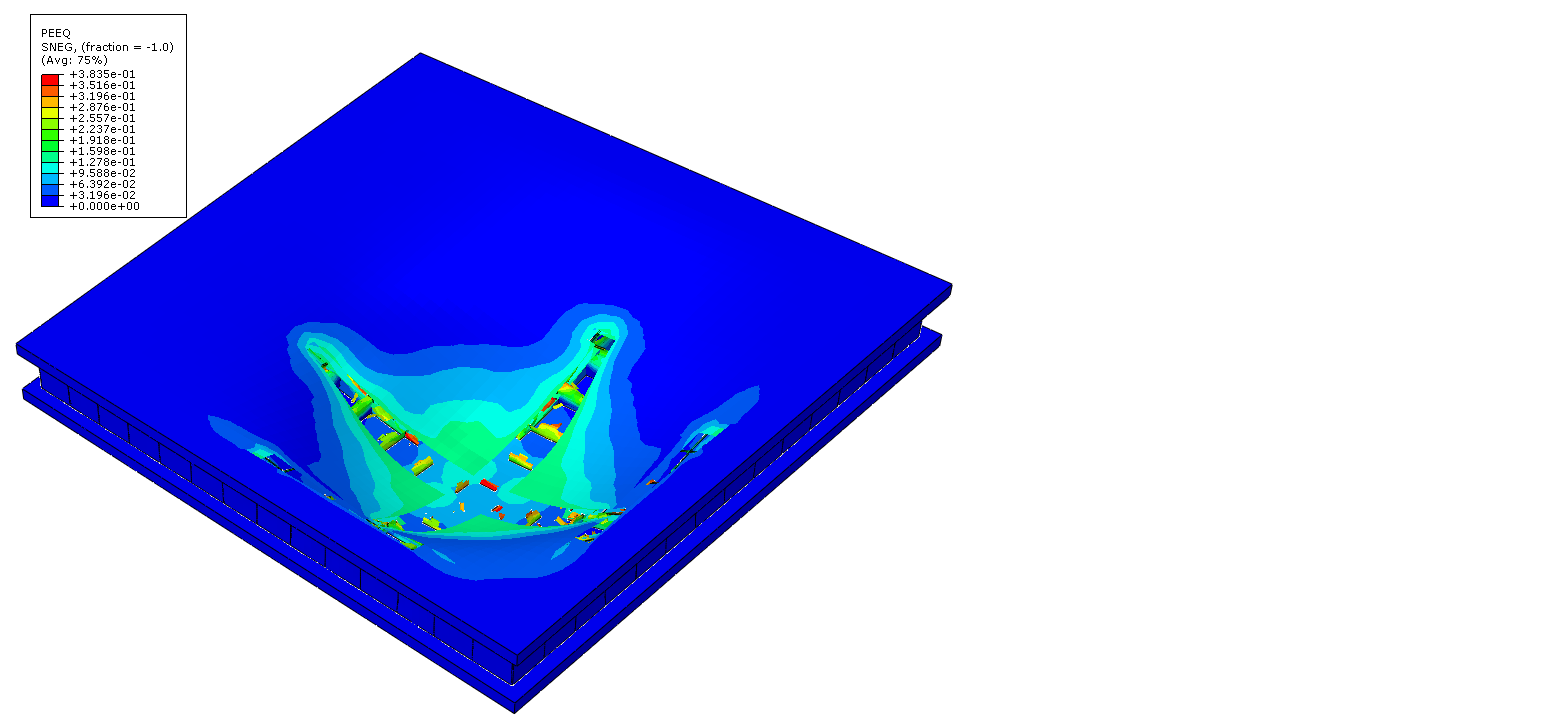}
         \label{fig:PEEQ}
     }
    \caption{FE model setup and deformation: \psubref{fig:BCs} Sandwich panel assembly with applied boundary conditions and loads. \psubref{fig:PEEQ} The deformed initial lattice design at the end of the simulation, colored by the equivalent plastic strain.}
    \label{cae}
\end{figure}

In the second example, a square lattice with side length of 100 mm and a height of 8.5 mm was considered. Two square face sheets with a side length of 105 mm and a thickness of 2.5 mm were placed on the top and bottom ends of the lattice, and perfect bonding was assumed between the lattice core and face sheets. The assembly was held fixed at the outer region defined by an offset width of 6.07 mm. Taking the center of the face sheet as the origin, 8 g of TNT was placed at position (15,15) mm with a stand-off distance of 19 mm from the top surface of the top face sheet. The blast loading generated by the TNT was modeled using the CONWEP model \cite{randers1997airblast}, and the simulation duration is 0.1 ms. The material of all components was Ti-6Al-4V. Due to the complex and dynamic nature of the applied loading, we used the rate-dependent Johnson-Cook plasticity and damage models to capture its material behavior. Key parameters for the model were adapted from the work of Wang and Shi \citep{wang2013validation} (see Table 2 therein), and are presented in \tref{JC_props}. \fref{fig:BCs} depicts the applied BCs and loads. 

\begin{table}[h!]
    \caption{Material parameters for John-Cook plasticity and damage models, Wang and Shi \cite{wang2013validation} }
    \small
    \centering
    \begin{tabular}{ccccccc}
    
    \multicolumn{7}{c}{Johnson-Cook plasticity model} \\
    \hline
    Name & \vline & \begin{tabular}{@{}c@{}}Yield\\stress\\(MPa)\end{tabular} & \begin{tabular}{@{}c@{}}Hardening\\coefficient\\(MPa)\end{tabular} &
    \begin{tabular}{@{}c@{}}Strain\\hardening\\exponent\end{tabular} &
    \begin{tabular}{@{}c@{}}Strain\\rate\\constant\end{tabular} &
    \begin{tabular}{@{}c@{}}Thermal\\softening\\exponent\end{tabular} \\
    \hline
    Value & \vline & 1098 & 1092 & 0.93 & 0.014 & 1.1 \\
    \hline
    
    \multicolumn{7}{c}{Johnson-Cook damage model} \\
    \hline
    Name & \vline & $d_1$ & $d_2$ & $d_3$ & $d_4$ & $d_5$ \\
    \hline
    Value & \vline & -0.09 & 0.27 & 0.48 & 0.014 & 3.87 \\
    \hline
    
    \end{tabular}
    \label{JC_props}
\end{table}

The lattice cross section was partitioned into 14 cells in the X- and Y-directions, and the out-most boundaries were nondesignable with a constant thickness of 0.5 mm. The initial thickness distribution is a uniform wall thickness of 0.25 mm. The optimization used the first thickness update scheme and maintained a mass identical to the initial design. \fref{fig:PEEQ} depicts the deformation and equivalent plastic strain distribution on the initial lattice design at the end of the simulation.

\begin{figure}[h!]
    \newcommand\x{0.3}
    \centering
    \begin{tabular}{ c  c  c  }
    \begin{minipage}[c]{\x\textwidth}
       \centering 
        \subfloat[Iteration 1]{\includegraphics[trim={2cm 0.49cm 2cm 1.5cm},clip,width=\textwidth]{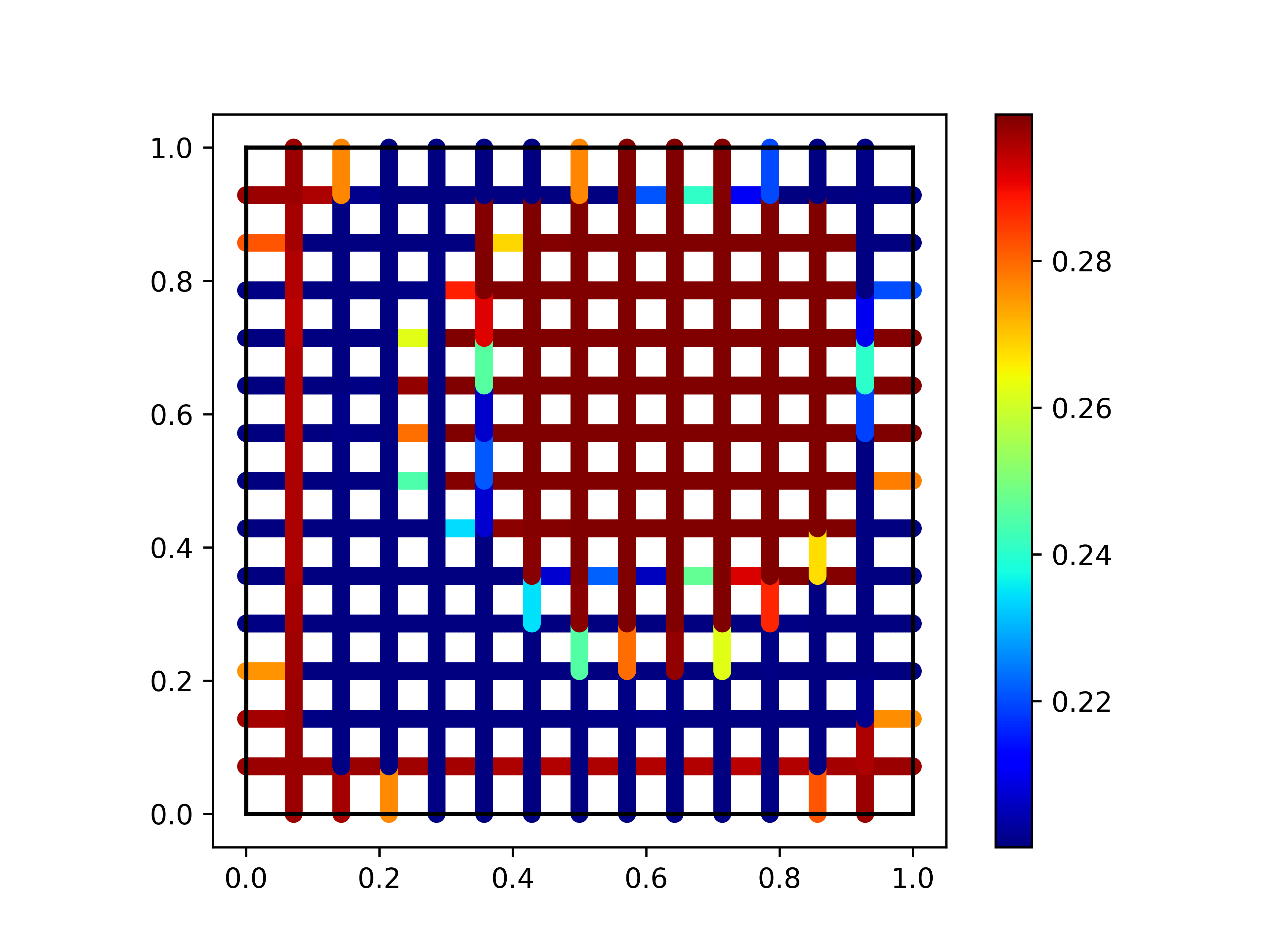}
        \label{fig:d1}}
    \end{minipage}
    &
    \begin{minipage}[c]{\x\textwidth}
       \centering 
        \subfloat[Iteration 3]{\includegraphics[trim={2cm 0.49cm 2cm 1.5cm},clip,width=\textwidth]{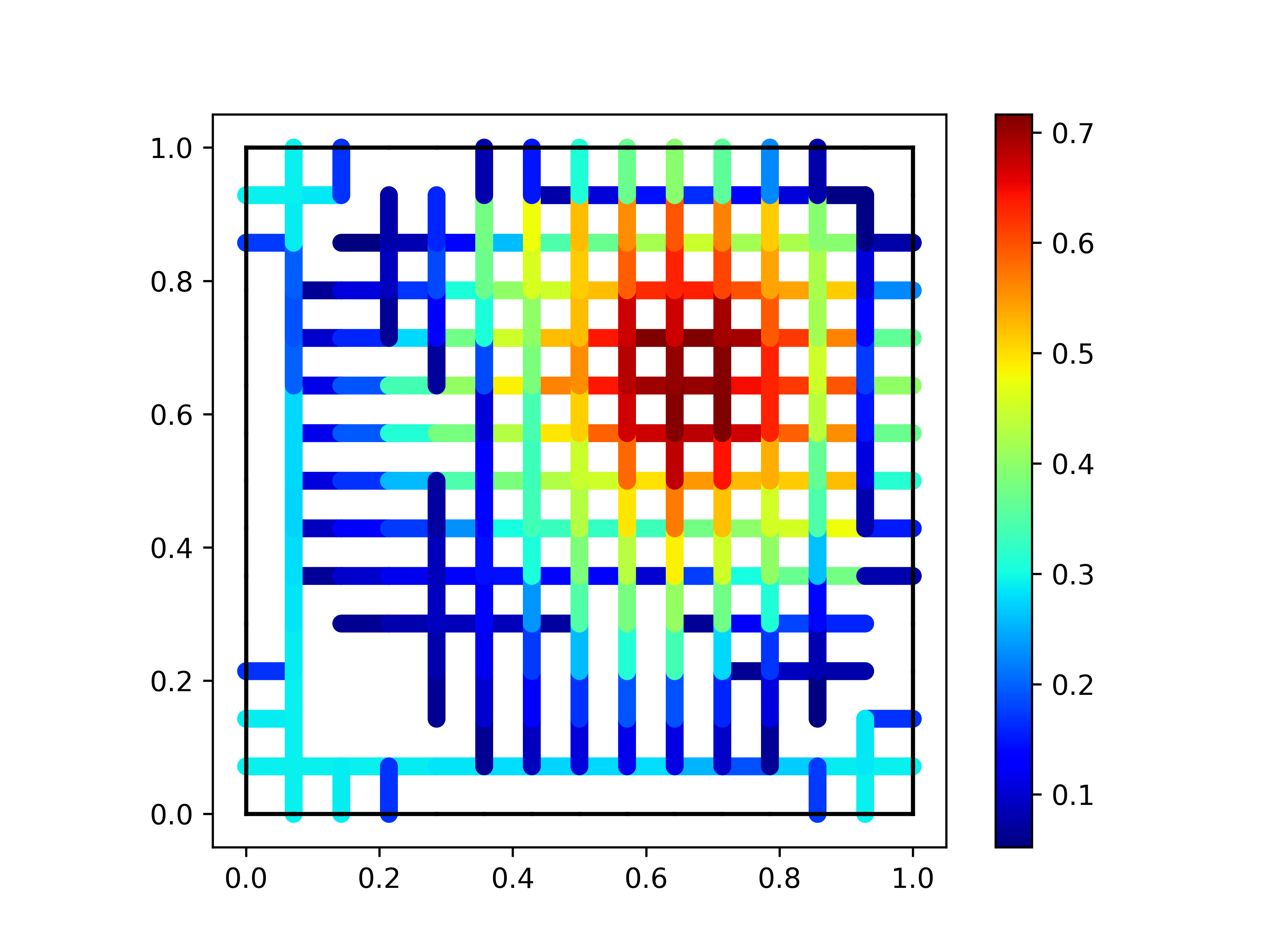}
        \label{fig:d2}}
    \end{minipage}\\

    \begin{minipage}[c]{\x\textwidth}
       \centering 
        \subfloat[Iteration 12]{\includegraphics[trim={2cm 0.49cm 2cm 1.5cm},clip,width=\textwidth]{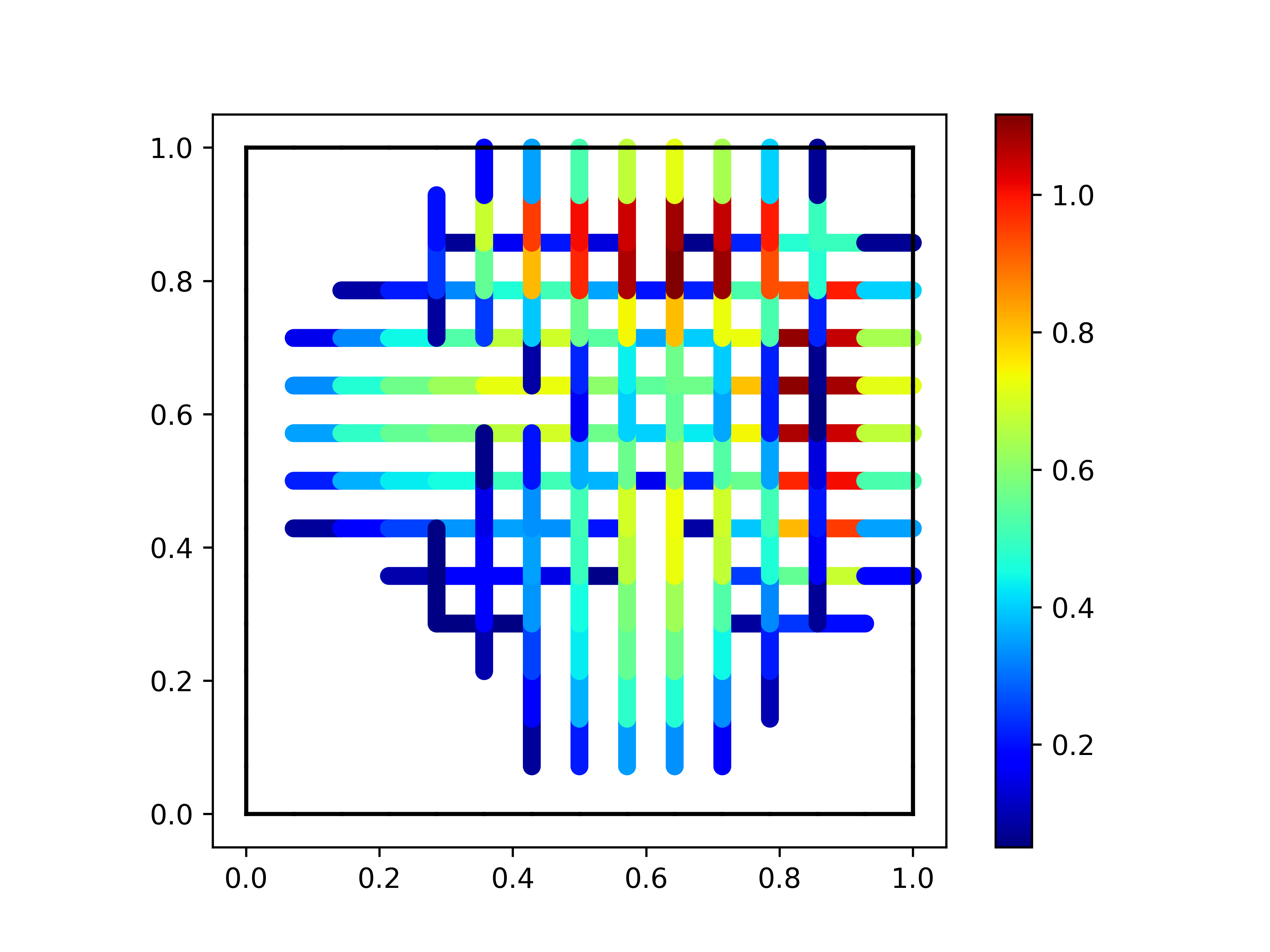}
        \label{fig:d3}}
    \end{minipage}
    &
    \begin{minipage}[c]{\x\textwidth}
       \centering 
        \subfloat[Iteration 24 $^*$]{\includegraphics[trim={2cm 0.49cm 2cm 1.5cm},clip,width=\textwidth]{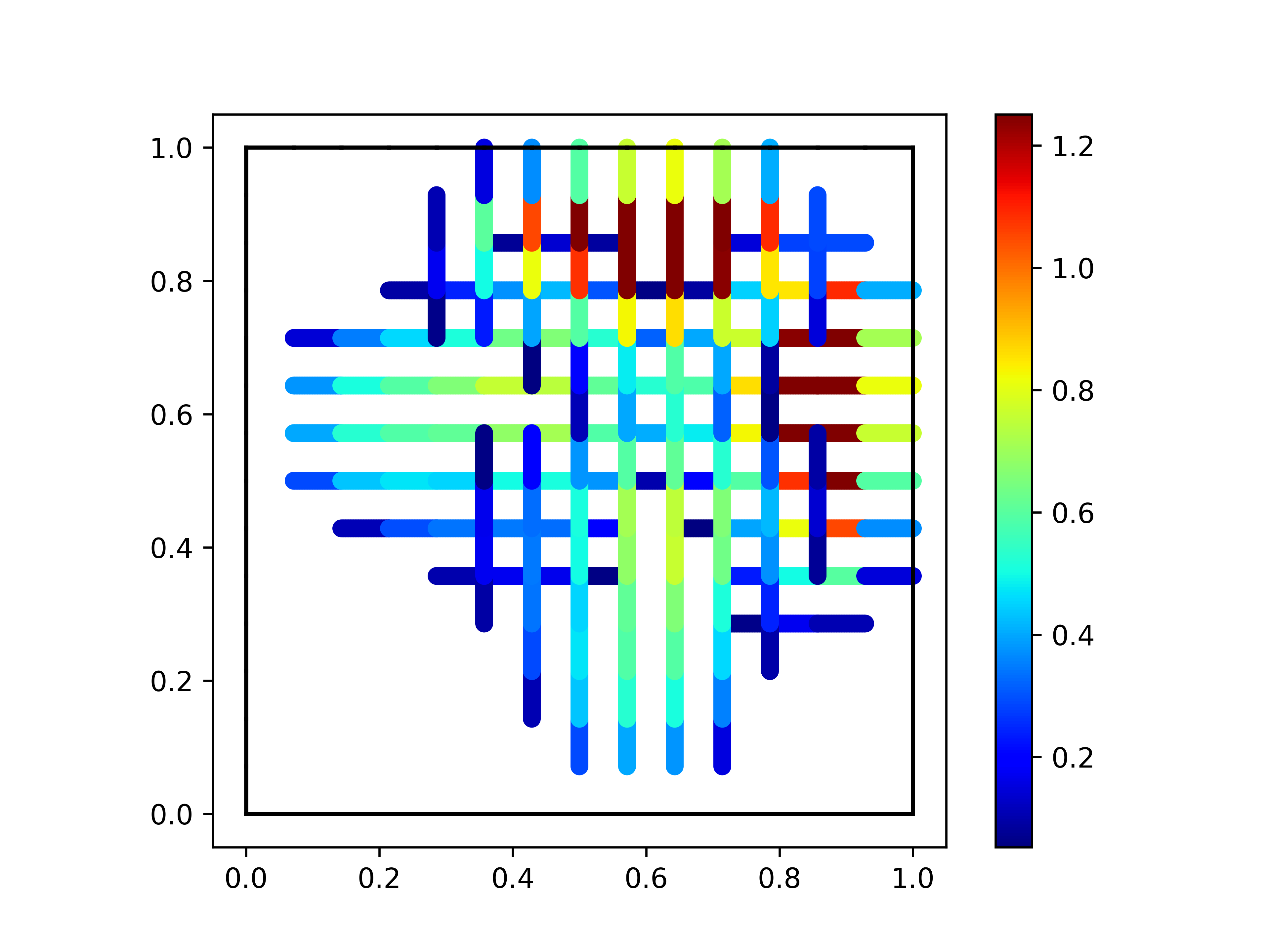}
        \label{fig:d4_opt}}
    \end{minipage}

    \end{tabular}
    \caption{Lattice core designs generated at different iterations, the optimized design that gives highest stiffness is marked with the superscript $^*$. Iteration 3 is when the update scheme first decided to remove walls. The color is assigned based on lattice wall thickness. The out-most boundary walls are marked as black lines. Note the similarity between iterations 12 and 24.}
    \label{case2_designs}
\end{figure}

The thickness update scheme did not converge on a design within the 0.01 mm tolerance, and ran the entire 25 design iterations. The intermediate designs generated by LatticeOPT are shown in \fref{case2_designs}. We clearly see that the update scheme decided to strengthen walls near the blast center, and quickly removed walls at the lower left corner, which were far away from the blast center. However, it is interesting to note that the top and right edges immediately connecting the lattice walls under the fixture were assigned the highest thickness (see \fref{fig:d4_opt}), which helps with connecting the center of the lattice structure to the surrounding fixtures to prevent shear failure at the edge of the fixture. We also noticed that the design topology stabilizes after about 13 iterations, and only the wall thicknesses fluctuated.

\begin{figure}[h!] 
    \centering
     \subfloat[]{
         \includegraphics[width=0.32\textwidth]{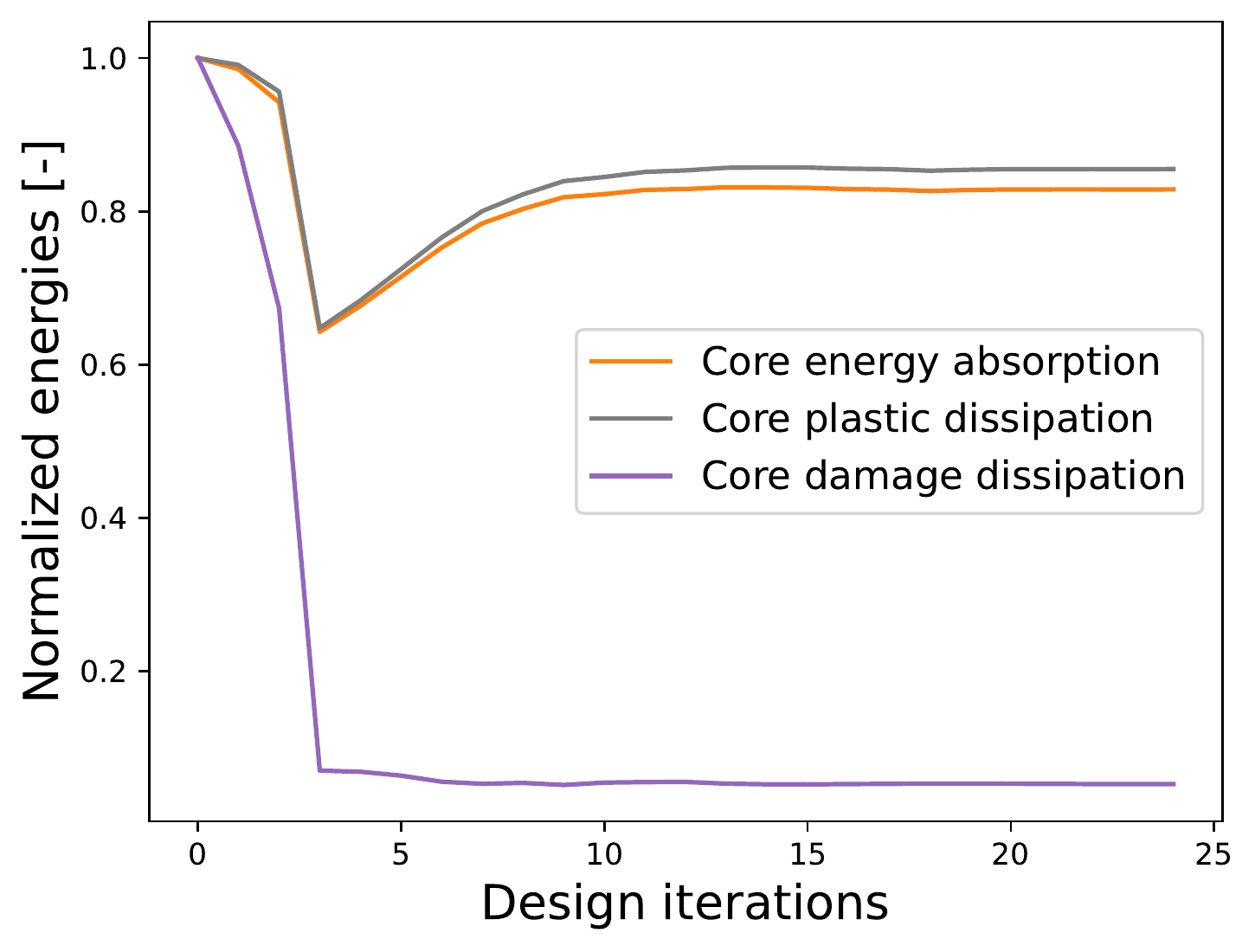}
         \label{fig:core_ene}
     }
     \subfloat[]{
         \includegraphics[width=0.32\textwidth]{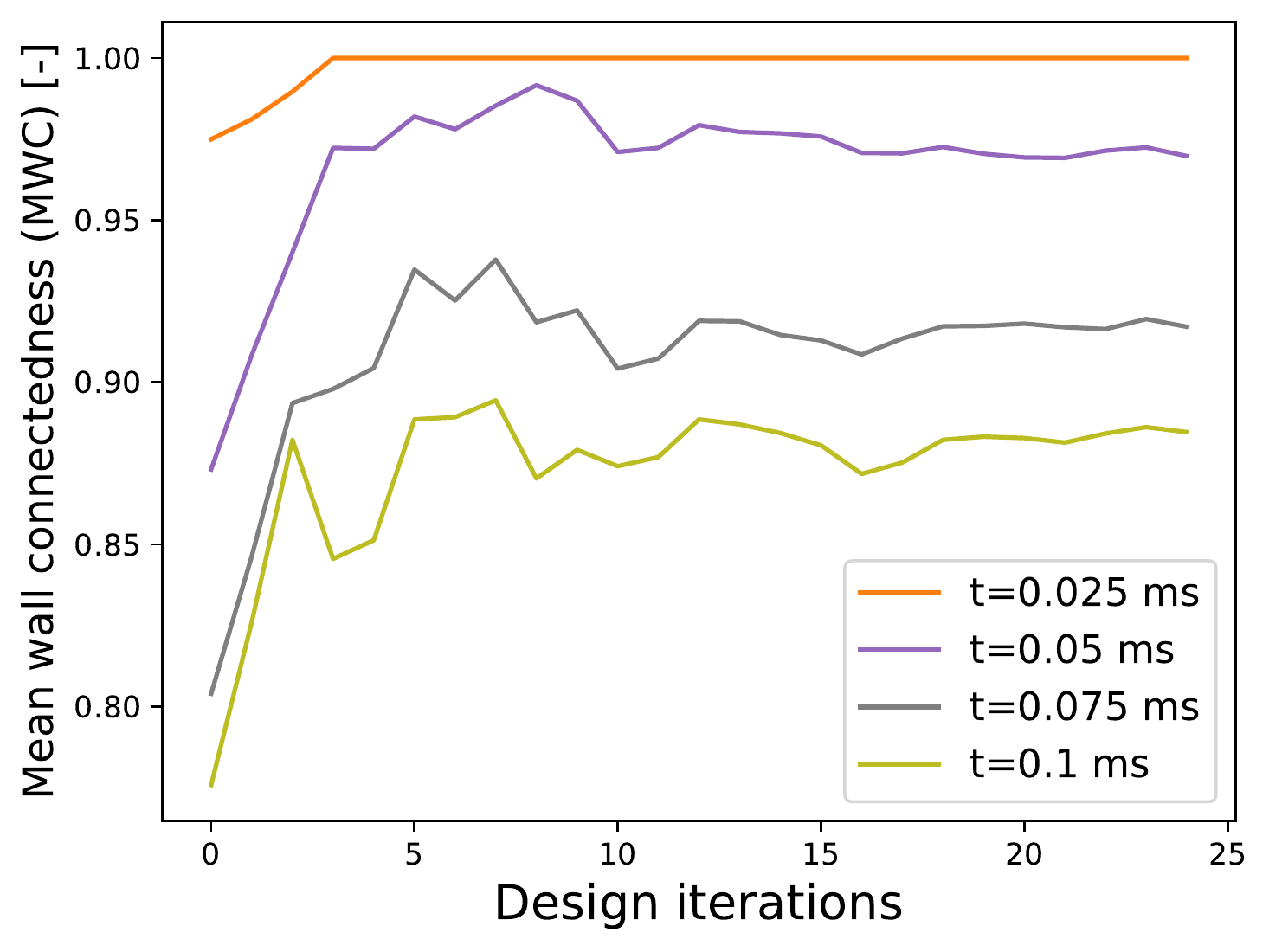}
         \label{fig:MWC}
     }
     \subfloat[]{
         \includegraphics[width=0.32\textwidth]{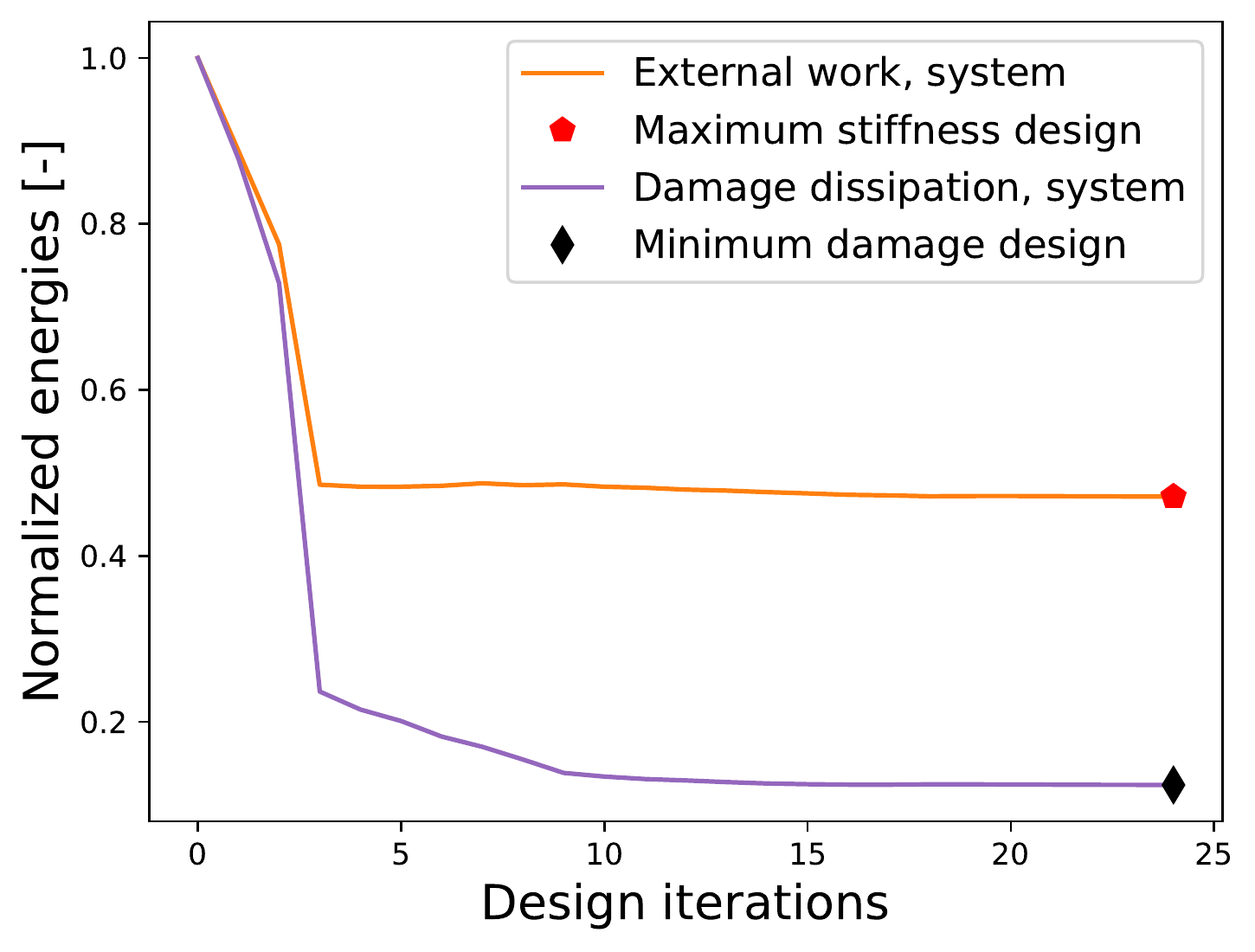}
         \label{fig:sys_ene}
     }
    \caption{Performance evaluation under blast loading: \psubref{fig:core_ene} Evolution history of normalized lattice core absorbed energy, plastic and damage dissipation. \psubref{fig:MWC} The mean wall connected metric evaluated at different time instances during the blast loading. \psubref{fig:sys_ene} Evolution history of normalized system absorbed energy, and damage dissipation. }
    \label{case2_core}
\end{figure}

For a critical evaluation of performance, we first investigate the lattice core, excluding the face sheets. For easiness of comparison, we normalized all energy quantities with the corresponding values of the initial design. \fref{case2_core} shows how the normalized energies and the MWC (defined in \eref{eq:MWC}) at four different time points evolved as the design optimization progressed. 
 
From \fref{fig:core_ene}, we see that as the design optimization progressed, the damage sustained by the lattice core was greatly reduced: the damage energy at iteration 24 was about 5.5\% of that in the initial design. Mitigation of core damage is also evident in \fref{fig:MWC}, where we see that the MWC at four different time points during the blast all showed an increasing trend with design iterations (a higher MWC indicates more walls remain interconnected). However, when inspecting the core energy absorption, we noticed that unlike the first example, the SEA at the end of the optimization is less than the initial design. This indicates that the apparent consequence of the thickness update schemes is not always to maximize structure SEA, and its implication will be clear once we inspect the sandwich panel system as a whole (including the face sheets).  

\begin{figure}[h!] 
    \centering
     \subfloat[]{
         \includegraphics[trim={0 0 18cm 0},clip,width=0.4\textwidth]{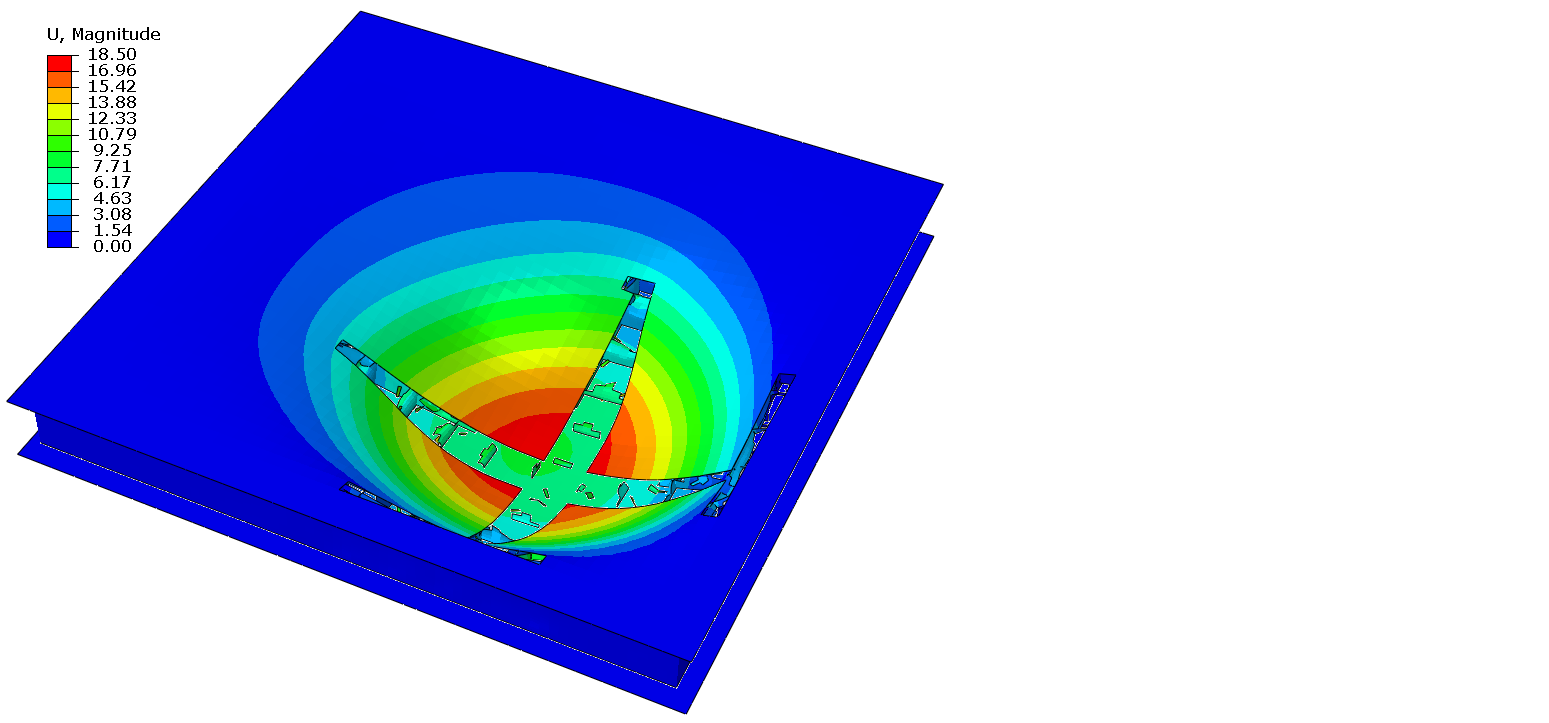}
         \label{fig:ini_def}
     }
     \subfloat[]{
         \includegraphics[trim={0 0 18cm 0},clip,width=0.4\textwidth]{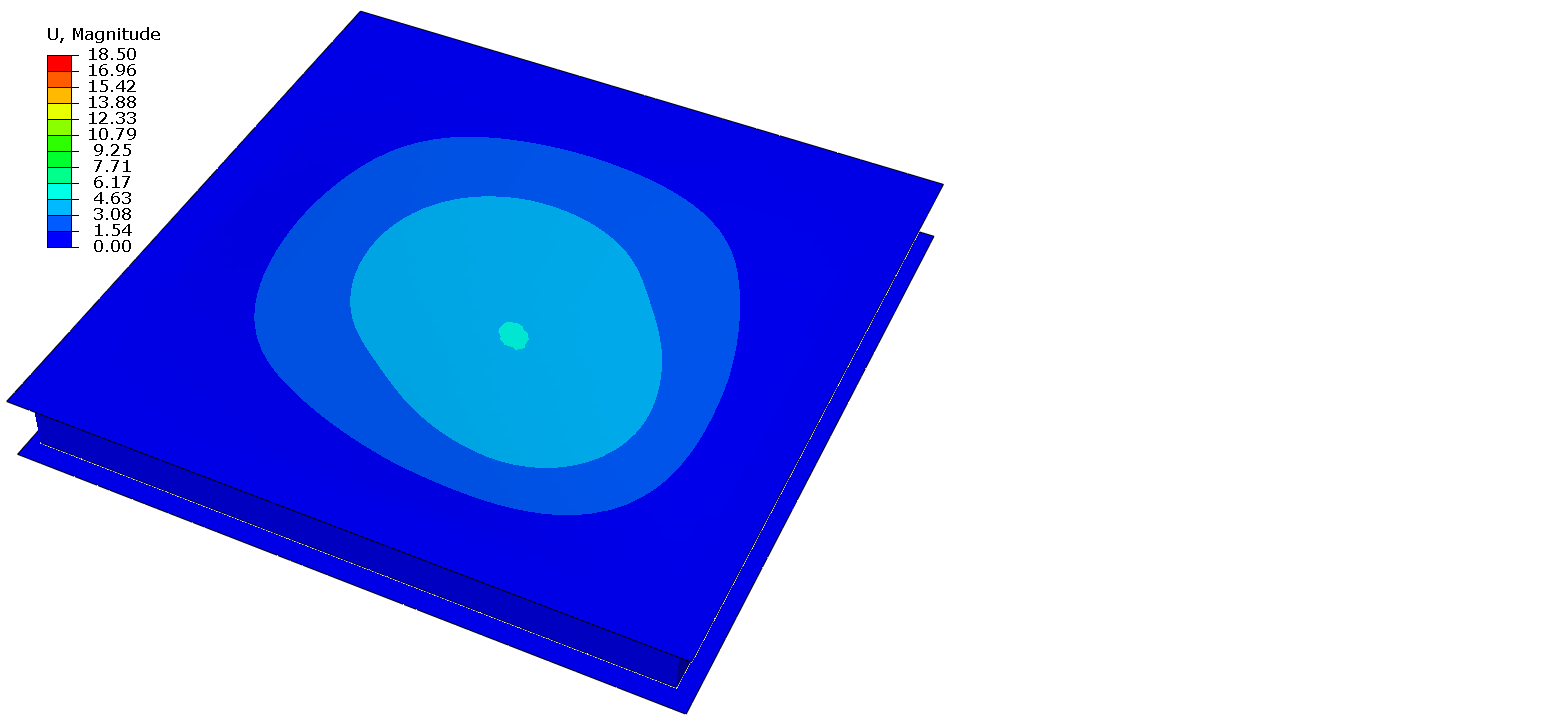}
         \label{fig:opt_def}
     }
    \caption{Deformation at the end of blast: \psubref{fig:ini_def} Initial lattice design. \psubref{fig:opt_def} Optimized lattice design (at design iteration 24).}
    \label{def_shape}
\end{figure}

For the sandwich panel system, we present the normalized total external work as well as damage dissipation in \fref{fig:sys_ene}, where it is evident that both the external work and damage dissipation decreased monotonically. The external work and damage dissipation at the last iteration was 47\% and 12.5\% of their initial values, respectively. Under a fixed load, decreasing external work is evidence that the system stiffness is increasing. Therefore, \fref{fig:sys_ene} elucidates that in the fixed-load scenario, the apparent consequence of the thickness update schemes is to increase system stiffness and minimize material damage. This example again highlights the effectiveness of LatticeOPT: in a case with complex blast loading and a total of 364 designable lattice walls, our framework was able to reduce the sandwich panel damage by 87.5\% in merely 25 design iterations. For a direct visual comparison of system stiffness and sustained damage, the deformed shape of the sandwich panel was plotted for the initial and optimized designs, shown in \fref{def_shape}.

\begin{figure}[h!] 
    \centering
     \subfloat[]{
         \includegraphics[width=0.32\textwidth]{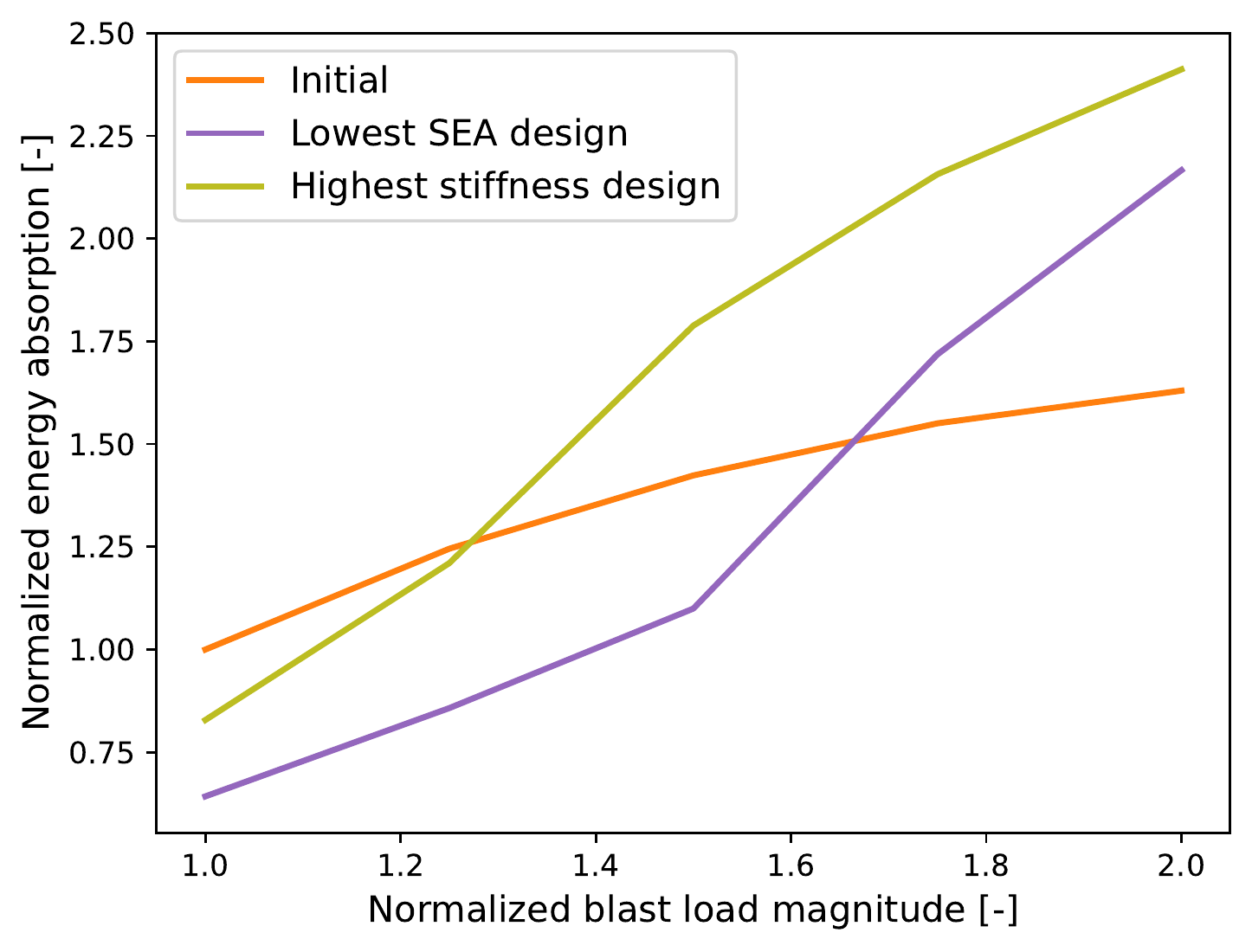}
         \label{fig:totE}
     }
     \subfloat[]{
         \includegraphics[width=0.32\textwidth]{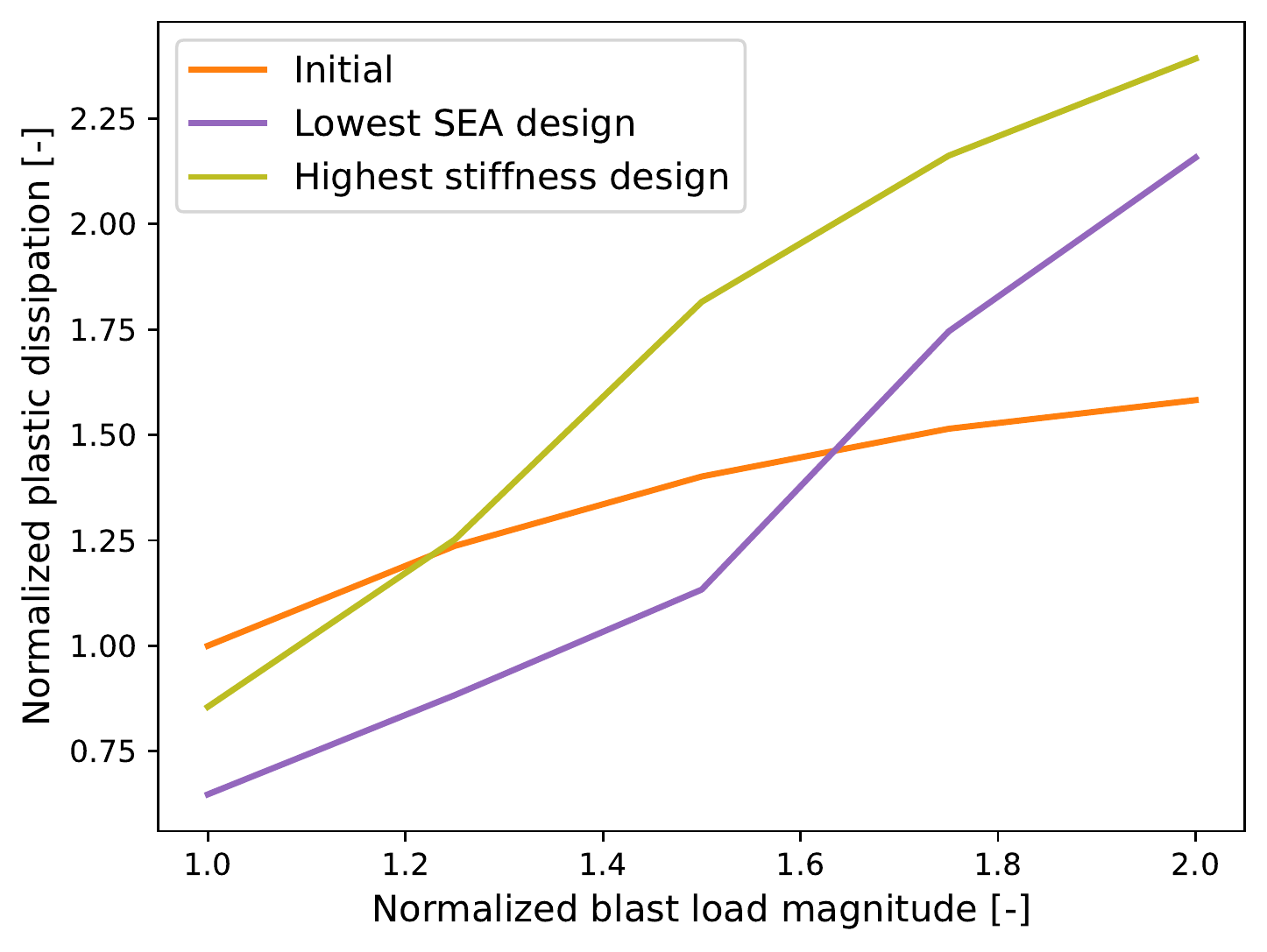}
         \label{fig:PD}
     }
     \subfloat[]{
         \includegraphics[width=0.32\textwidth]{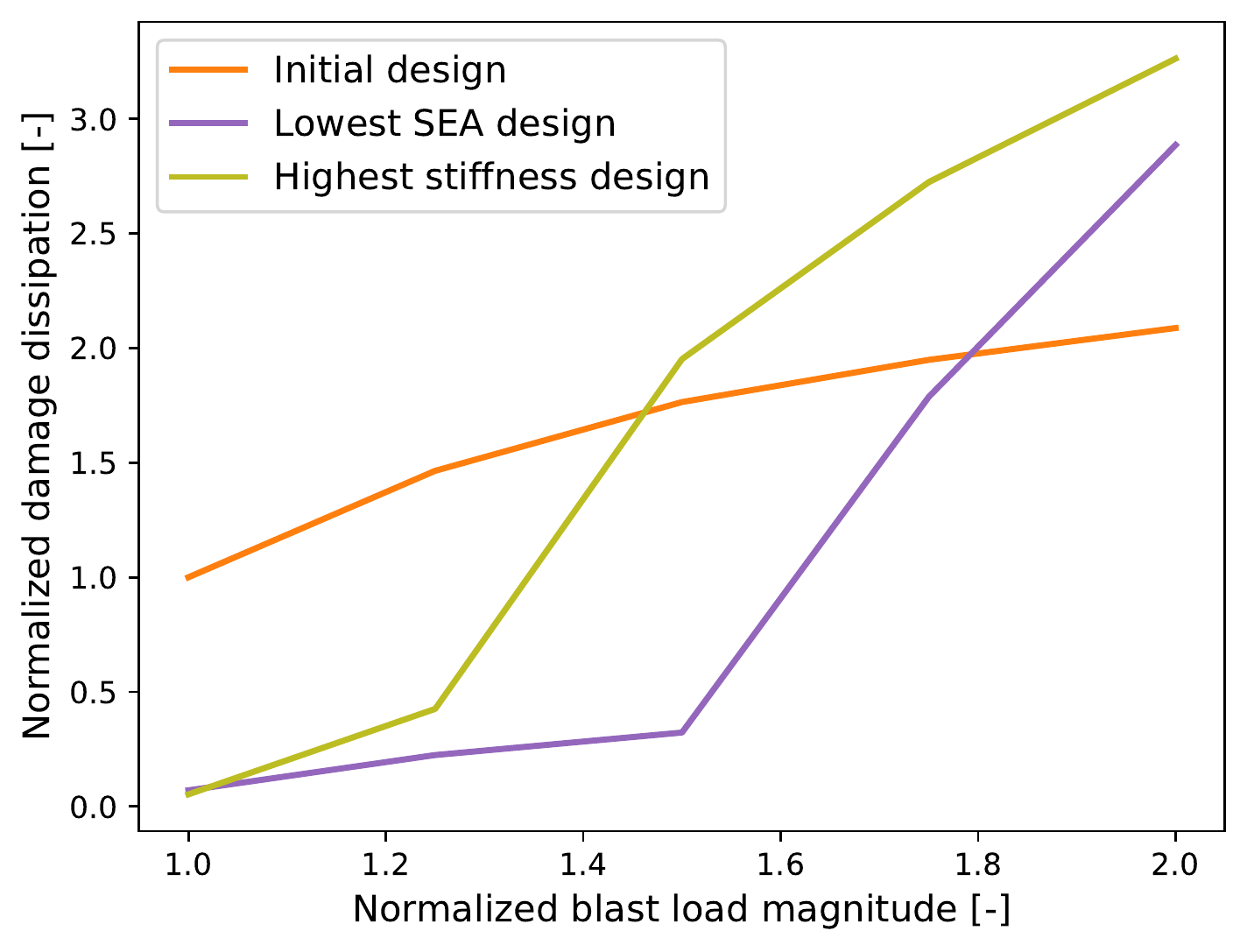}
         \label{fig:DMD}
     }
    \caption{Effects of increasing load magnitude: \psubref{fig:totE} Lattice core energy absorption. \psubref{fig:PD} Lattice core plastic dissipation. \psubref{fig:DMD} Lattice core damage dissipation. All quantities are normalized by that of the initial design at the original load level.}
    \label{increase_load}
\end{figure}

Finally, it is of practical interest to investigate how different designs respond as the magnitude of the blast load increases. Therefore, we gradually increased the normalized load magnitude from 1 to 2, and tested the initial lattice design, as well as design iterations 3 (has lowest SEA) and 24 (has highest stiffness), the plot of the normalized energies is presented in \fref{increase_load}.
We clearly see that, although the optimized design and the minimum SEA design had low energy absorption at the original load level, they were able to absorb much more energy than the initial design at elevated load magnitudes. At twice the original load magnitude, the optimized design and the minimum SEA design absorbed 48.3\% and 32.8\% more energy than the initial design. From \fref{fig:PD} and \fref{fig:DMD}, it is evident that the increased energy absorption is due to increases in plastic and damage dissipation for the optimized designs. This demonstrates that, besides stiffness enhancement, the optimized designs generated by LatticeOPT are also more robust, meaning that they have potential to absorb more energy should the load magnitude increase.

\subsection{Generation of periodic lattices}
\label{sec:ex3}
\begin{figure}[h!] 
    \centering
     \subfloat[]{
         \includegraphics[width=0.32\textwidth]{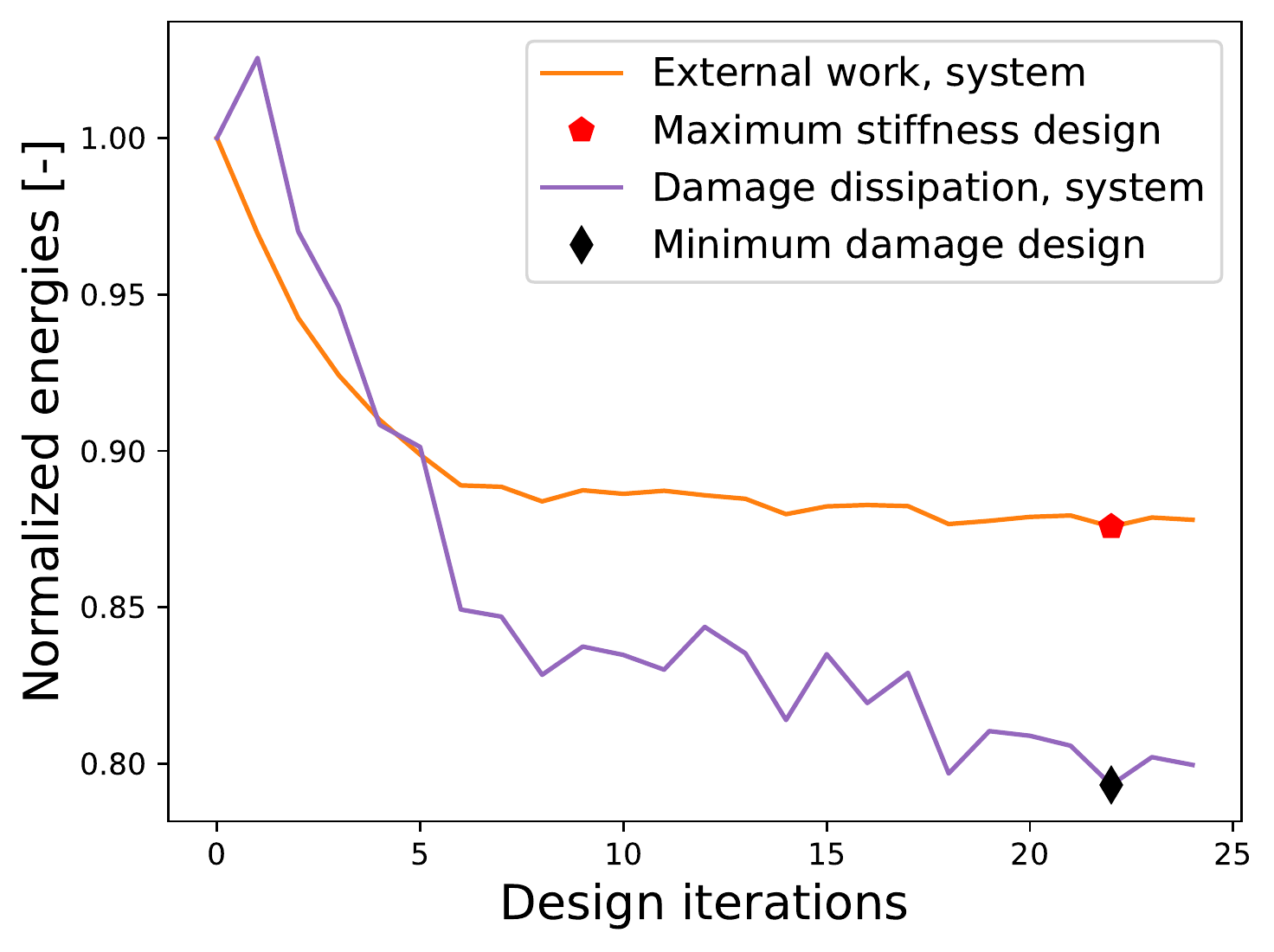}
         \label{fig:system_e_p}
     }
     \subfloat[Iteration 1]{
         \includegraphics[trim={2cm 0.49cm 2cm 1.5cm},clip,width=0.32\textwidth]{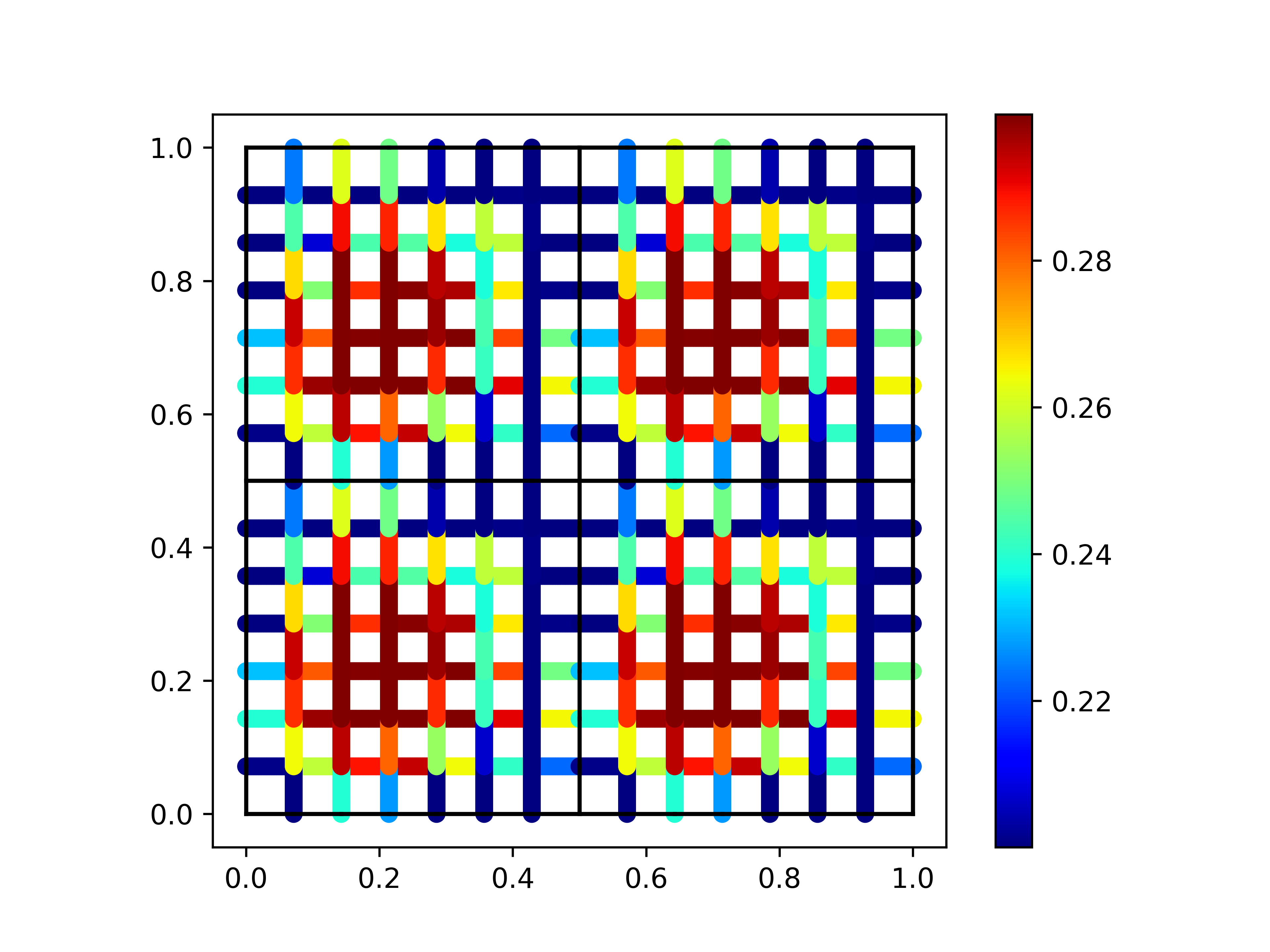}
         \label{fig:itr1}
     }
     \subfloat[Iteration 22 $^*$]{
         \includegraphics[trim={2cm 0.49cm 2cm 1.5cm},clip,width=0.32\textwidth]{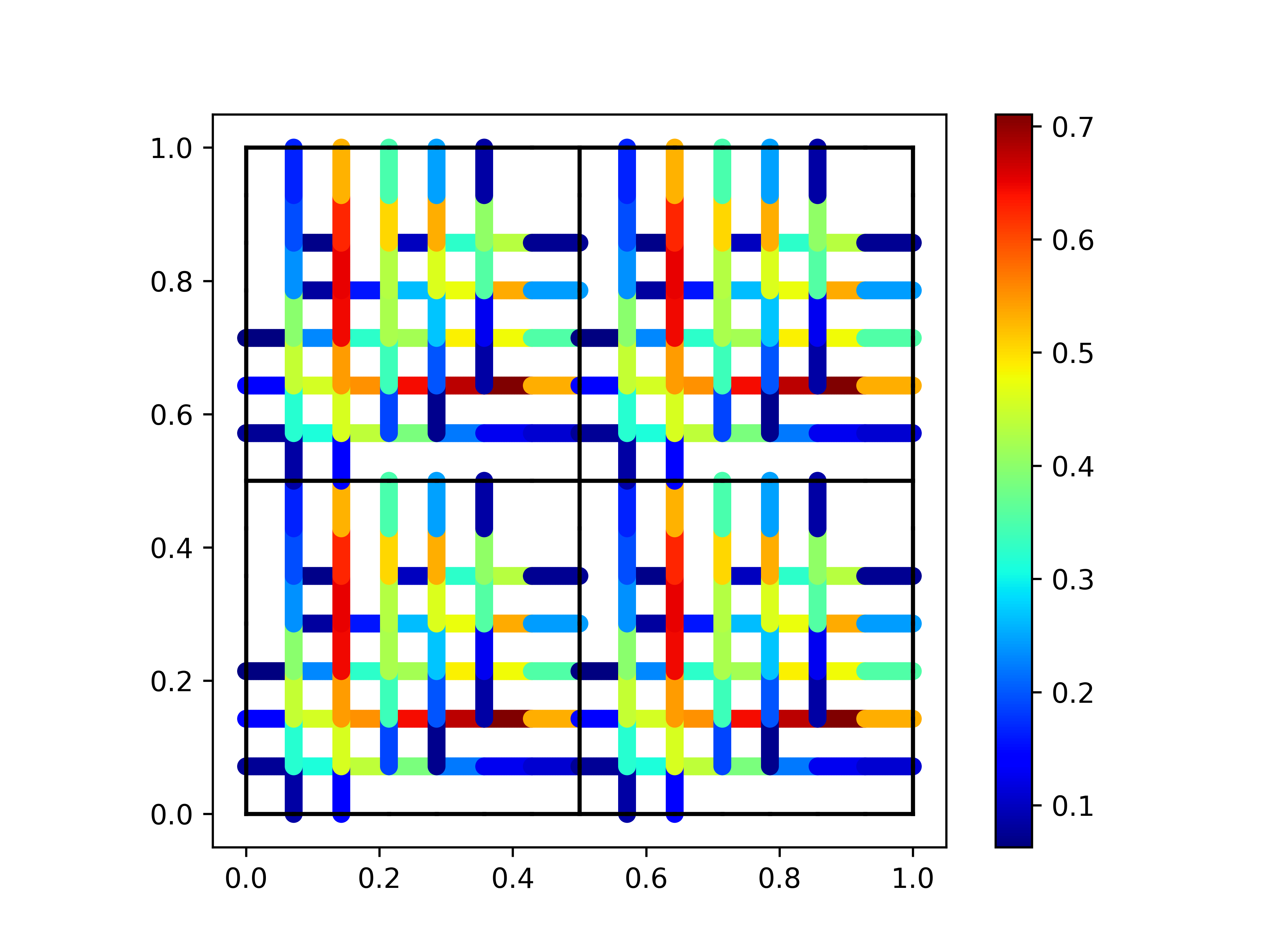}
         \label{fig:itr22}
     }
    \caption{Performance and designs of the periodic lattice core: \psubref{fig:system_e_p} External work and damage dissipation of the system. \psubref{fig:itr1} Design iteration 1. \psubref{fig:itr22} Design iteration 22, this design produced the highest system stiffness. Note the clear periodicity of the generated designs.}
    \label{periodic}
\end{figure}

The setup of the third example followed exactly from \sref{sec:ex2}, except we enforced a 2-by-2 periodic cell arrangement. The total external work and the damage dissipation of the sandwich panel assembly are shown in \fref{fig:system_e_p}, and the designs at iterations 1 and 22 (highest stiffness) are shown in \fref{fig:itr1} and \fref{fig:itr22}.

From \fref{fig:system_e_p}, we again confirmed that the apparent consequence of the thickness update scheme in this fixed-load case is to maximize system stiffness and minimize damage. However, we noticed that the periodic design, having a minimum damage of 79.3\% of the initial value, is inferior to the free design presented in \sref{sec:ex2}, which has a minimum damage of only 5.5\% of initial value. This is reasonable, as the 2-by-2 periodicity constraint greatly reduced the size of the design space, and a more restrained system typically offers inferior performance compared to a system with more design degrees of freedom (i.e., number of lattice walls). Inspecting the optimized designs, we clearly see the 2-by-2 periodic unit cell arrangement. Note the generated designs are distinct from those in \sref{sec:ex2} due to this added constraint.

\begin{figure}[h!] 
    \centering
     \subfloat[]{
         \includegraphics[width=0.32\textwidth]{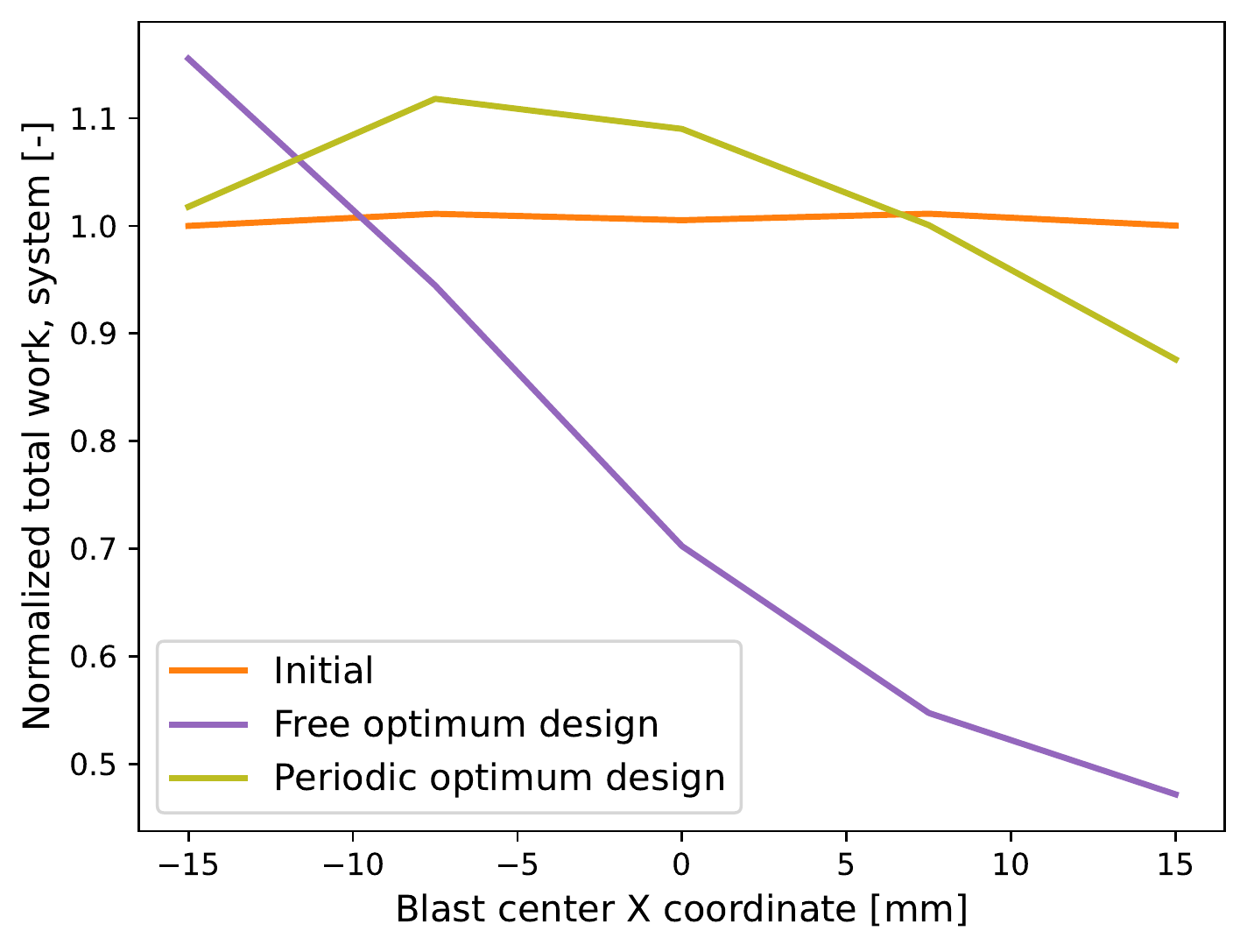}
         \label{fig:totWk}
     }
     \subfloat[]{
         \includegraphics[width=0.32\textwidth]{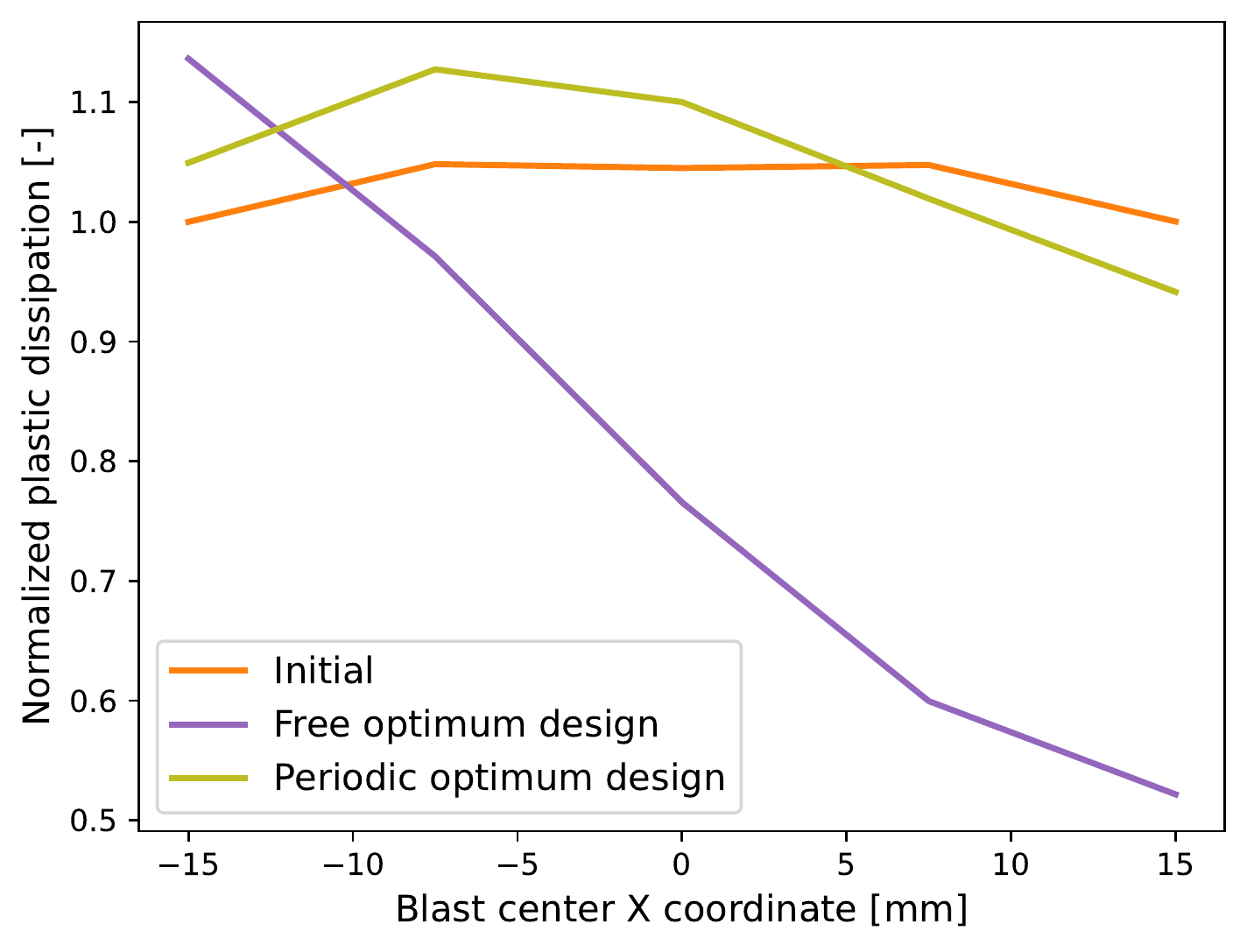}
         \label{fig:sPD}
     }
     \subfloat[]{
         \includegraphics[width=0.32\textwidth]{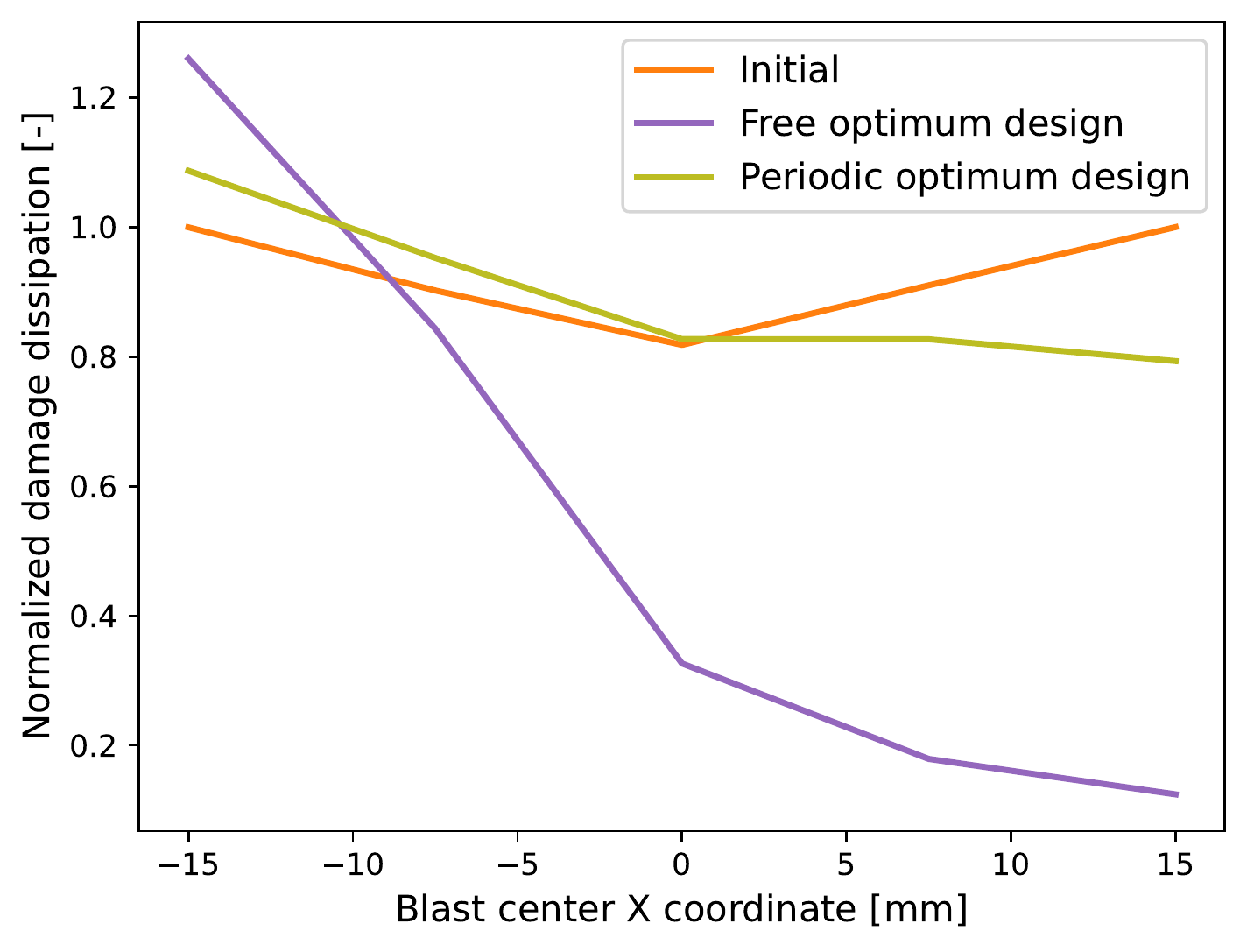}
         \label{fig:sDMD}
     }
    \caption{Effects of shifting blast center location: \psubref{fig:totWk} External work done to the sandwich panel. \psubref{fig:sPD} Sandwich panel plastic dissipation. \psubref{fig:sDMD} Sandwich panel damage dissipation. All quantities are normalized by that of the initial design at the original blast center location. The original blast center location (15,15) mm is shown on the right of each plot.}
    \label{shifted_loc}
\end{figure}

\begin{figure}[h!]
    \newcommand\x{0.28}
    \centering
    \begin{tabular}{ c  c  c  }
    \begin{minipage}[c]{\x\textwidth}
       \centering 
        \subfloat[Initial, $X_{bc}$ = -15 mm]{\includegraphics[trim={10cm 0 14cm 0},clip,width=\textwidth]{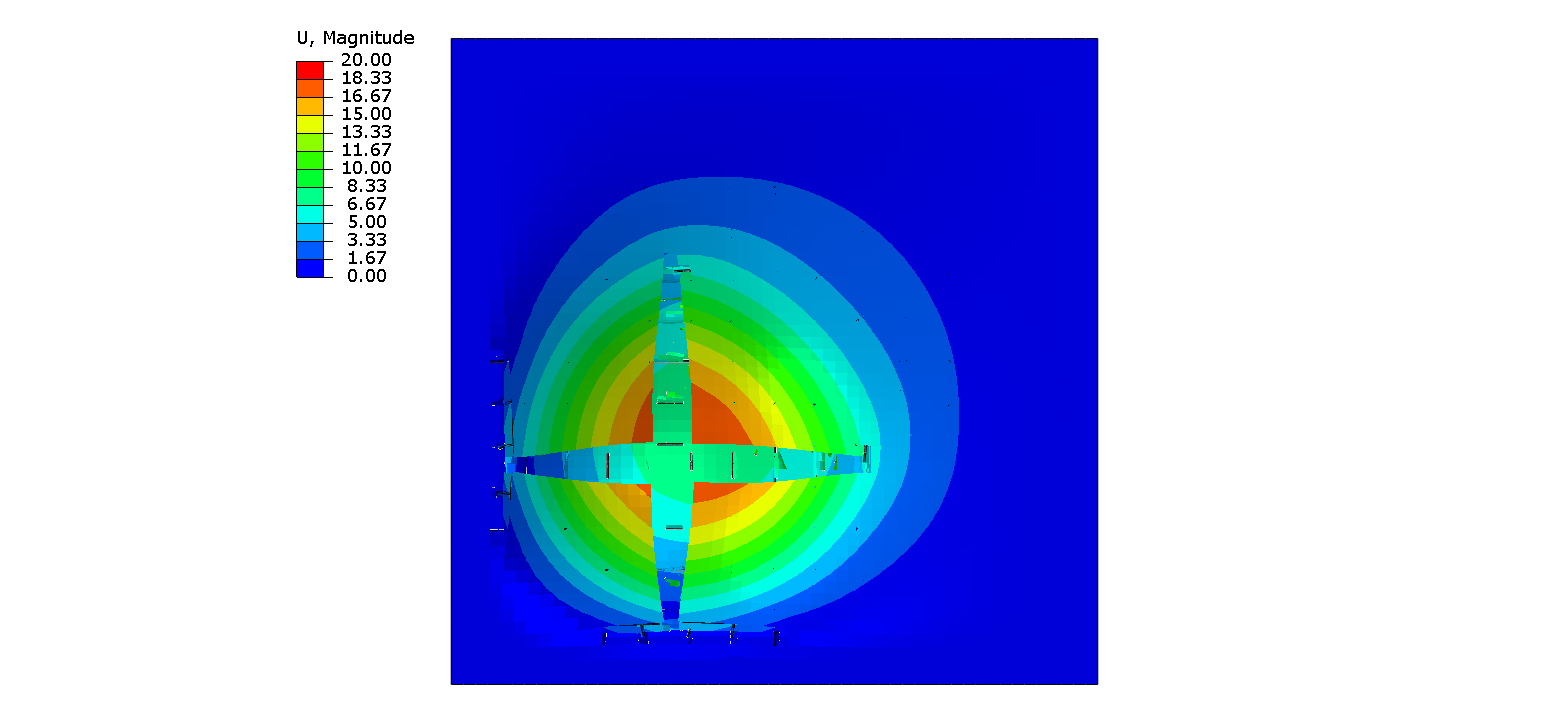}
        \label{fig:d1}}
    \end{minipage}
    &
    \begin{minipage}[c]{\x\textwidth}
       \centering 
        \subfloat[Initial, $X_{bc}$ = 0 mm]{\includegraphics[trim={10cm 0 14cm 0},clip,width=\textwidth]{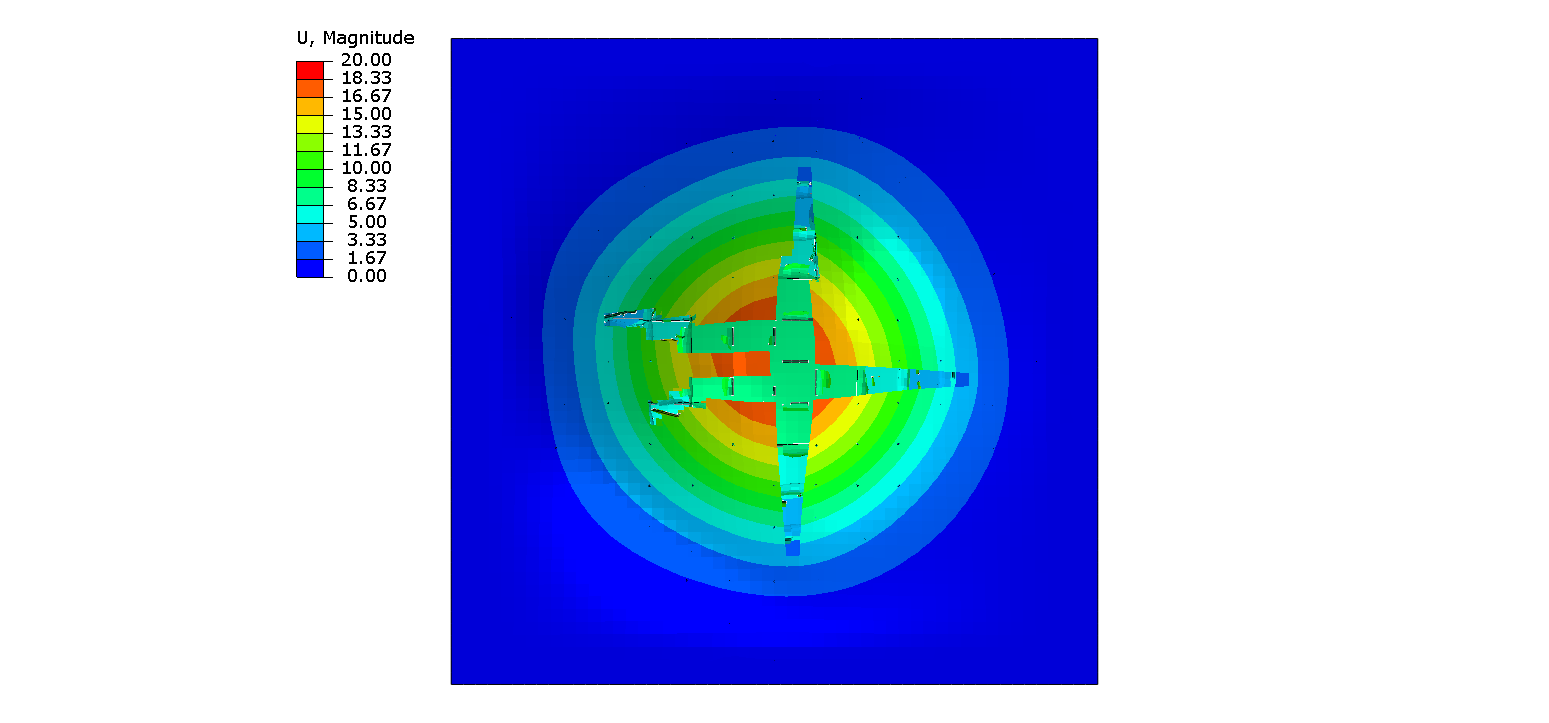}
        \label{fig:d2}}
    \end{minipage}
    &
    \begin{minipage}[c]{\x\textwidth}
       \centering 
        \subfloat[Initial, $X_{bc}$ = 15 mm]{\includegraphics[trim={10cm 0 14cm 0},clip,width=\textwidth]{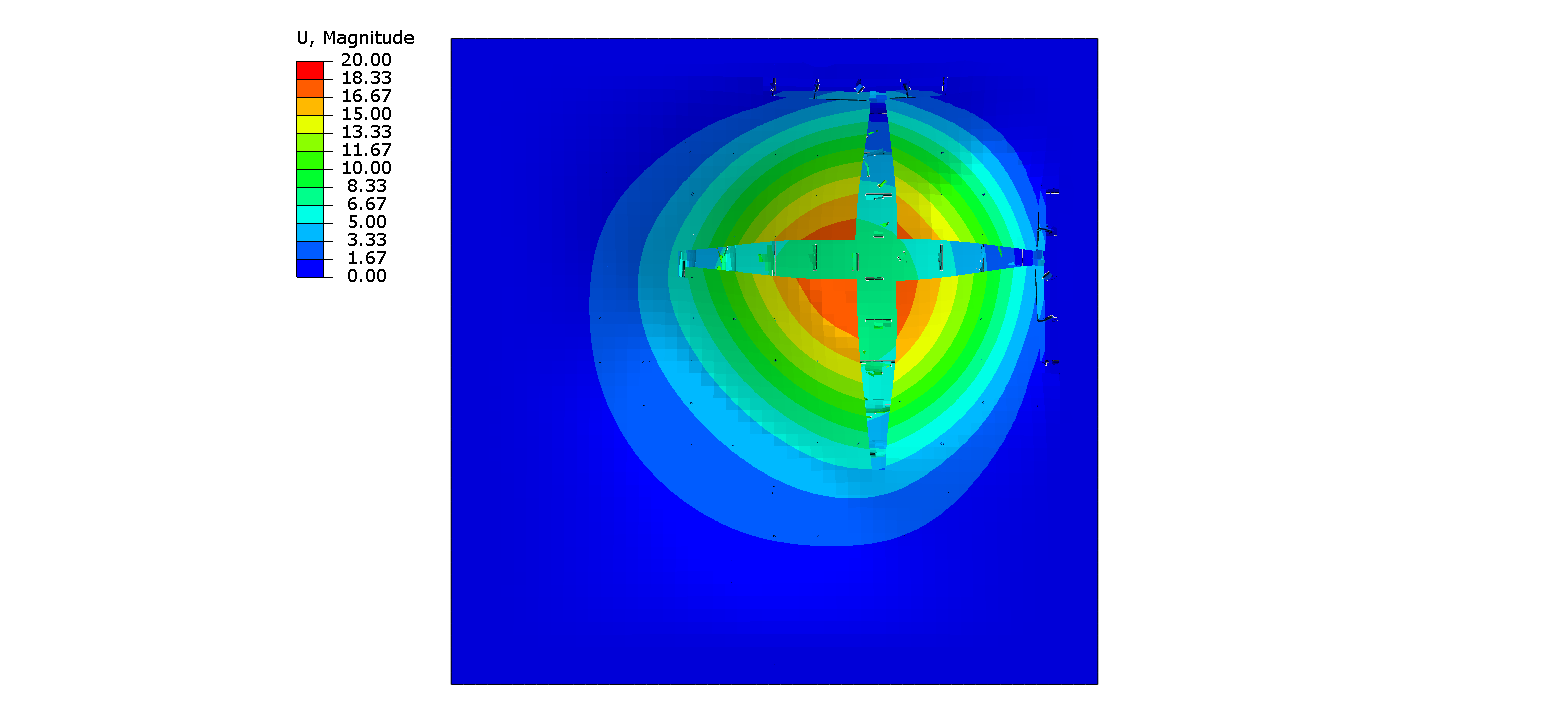}
        \label{fig:d3}}
    \end{minipage}\\

    \begin{minipage}[c]{\x\textwidth}
       \centering 
        \subfloat[Free, $X_{bc}$ = -15 mm]{\includegraphics[trim={10cm 0 14cm 0},clip,width=\textwidth]{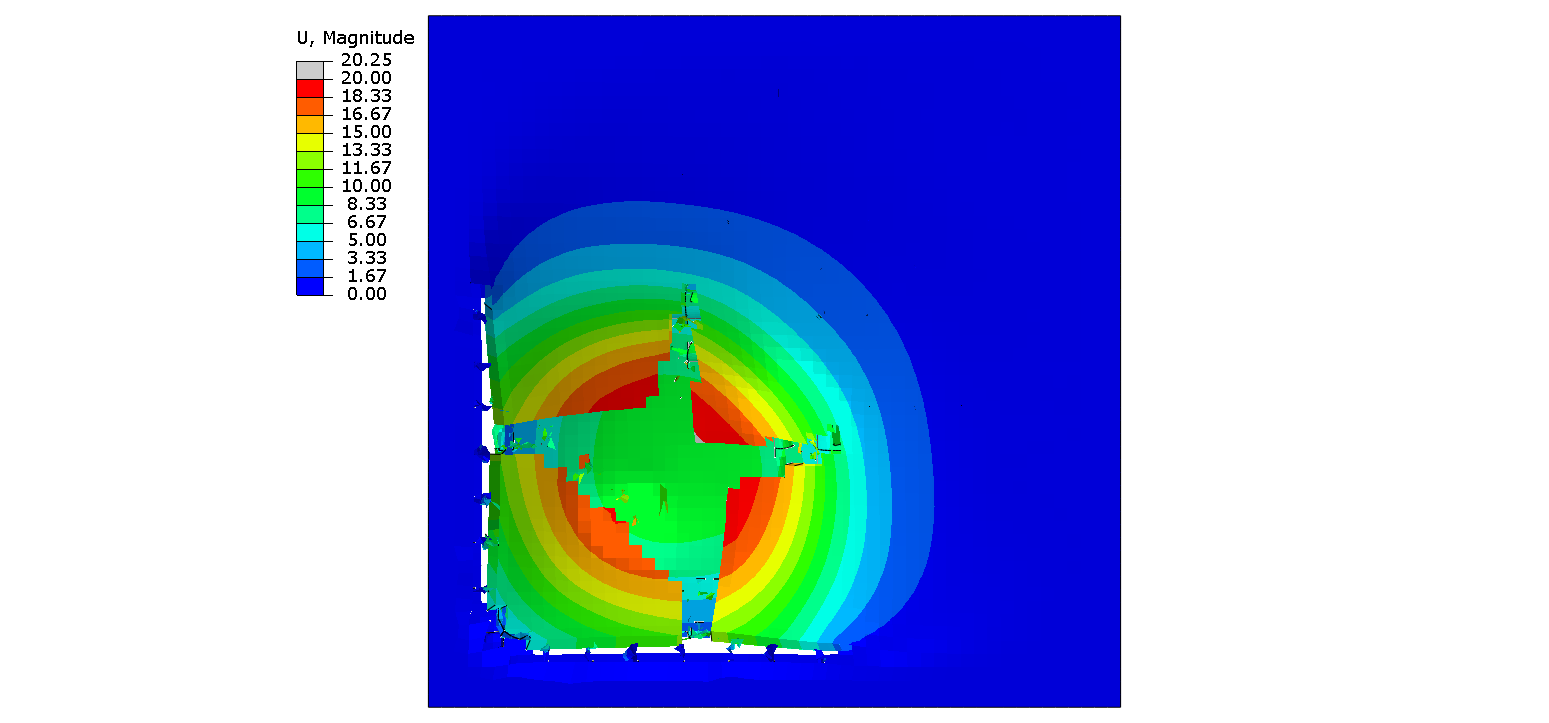}
        \label{fig:d4}}
    \end{minipage}
    &
    \begin{minipage}[c]{\x\textwidth}
       \centering 
        \subfloat[Free, $X_{bc}$ = 0 mm]{\includegraphics[trim={10cm 0 14cm 0},clip,width=\textwidth]{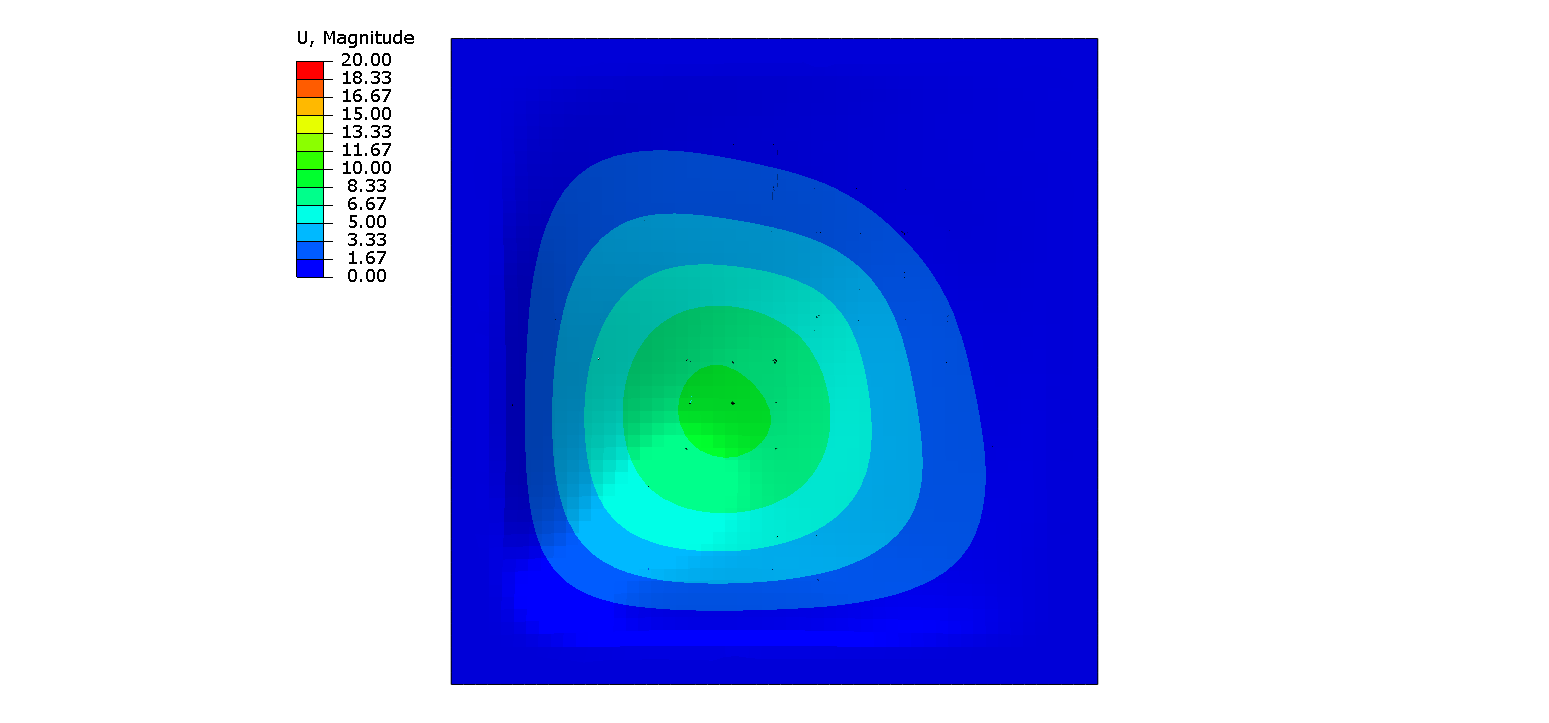}
        \label{fig:d5}}
    \end{minipage}
    &
    \begin{minipage}[c]{\x\textwidth}
       \centering 
        \subfloat[Free, $X_{bc}$ = 15 mm]{\includegraphics[trim={10cm 0 14cm 0},clip,width=\textwidth]{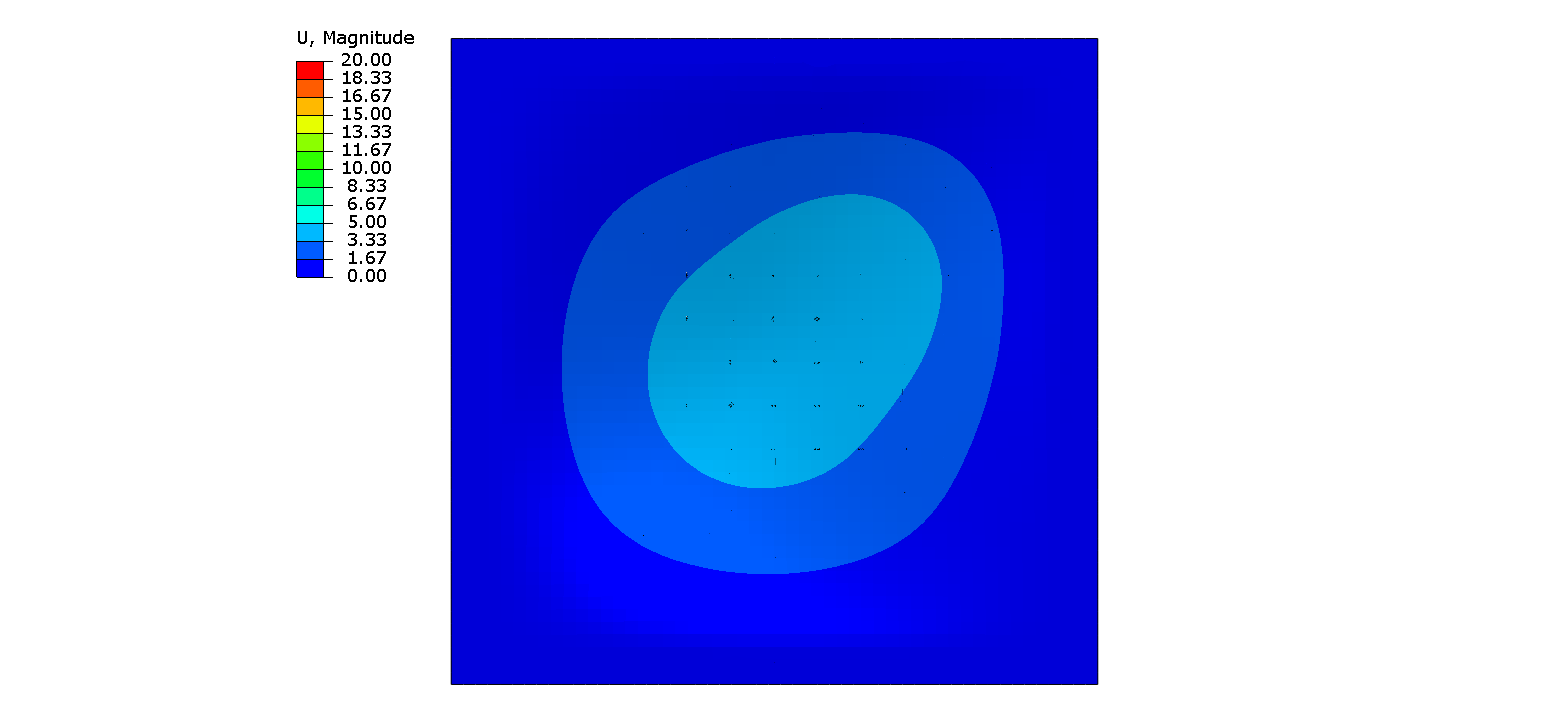}
        \label{fig:d6}}
    \end{minipage}\\
    
    \begin{minipage}[c]{\x\textwidth}
       \centering 
        \subfloat[Periodic, $X_{bc}$ = -15 mm]{\includegraphics[trim={10cm 0 14cm 0},clip,width=\textwidth]{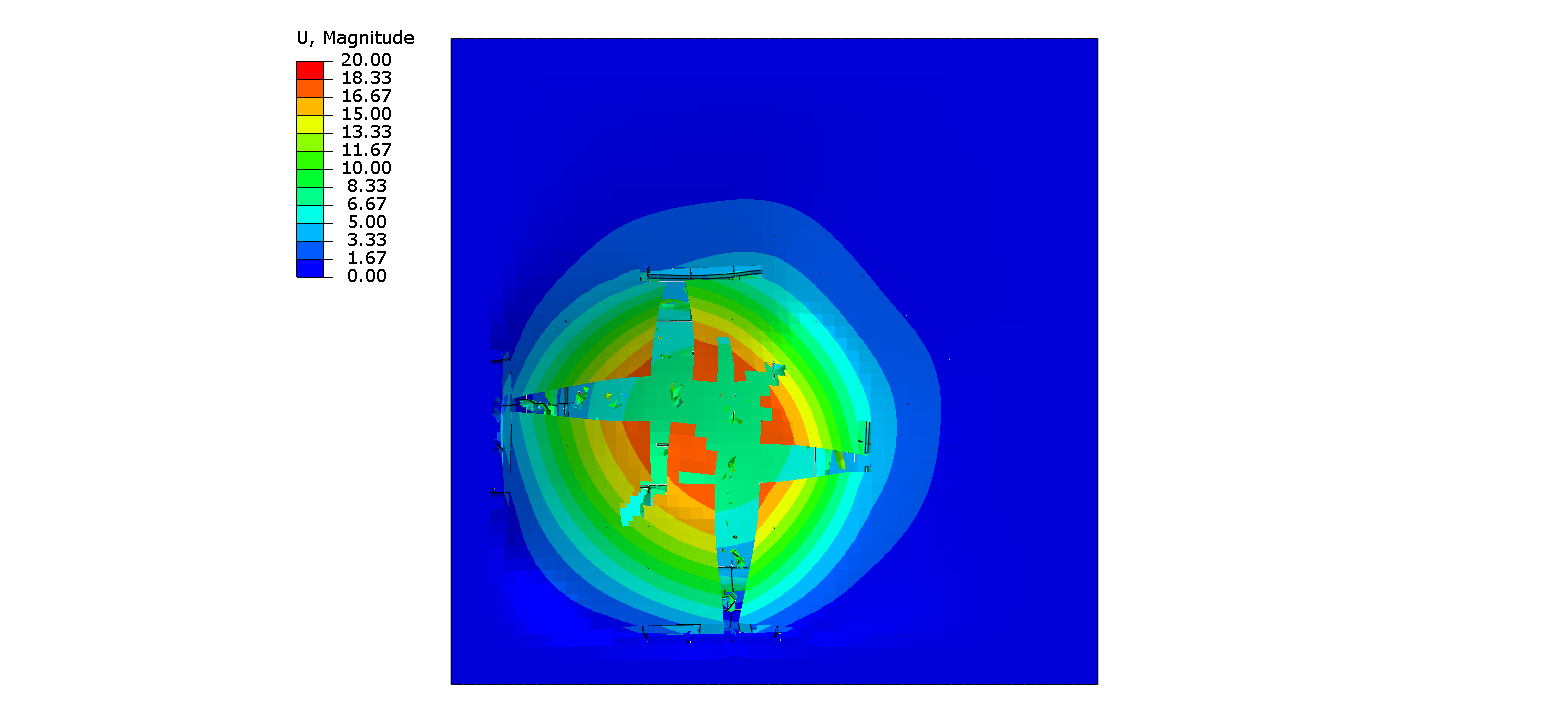}
        \label{fig:d7}}
    \end{minipage}
    &
    \begin{minipage}[c]{\x\textwidth}
       \centering 
        \subfloat[Periodic, $X_{bc}$ = 0 mm]{\includegraphics[trim={10cm 0 14cm 0},clip,width=\textwidth]{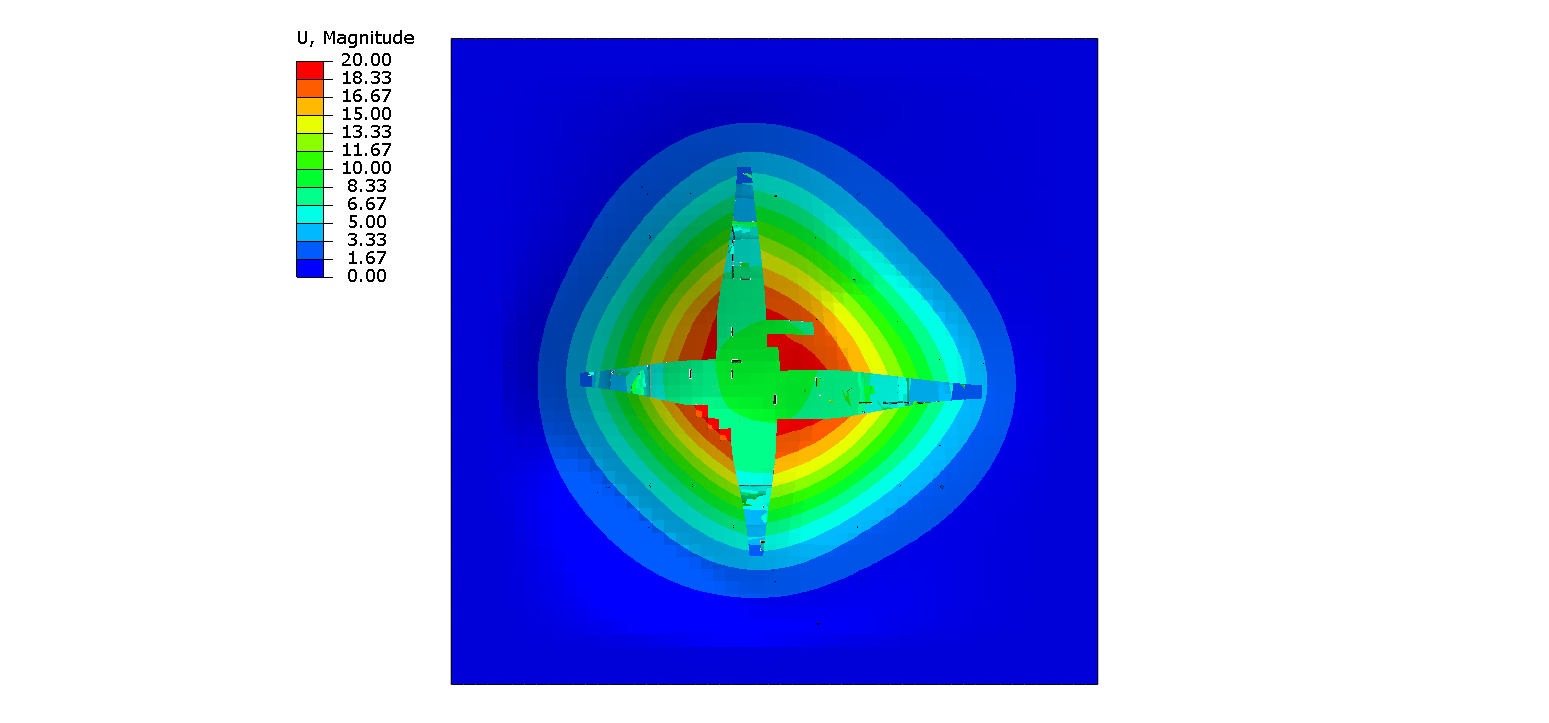}
        \label{fig:d8}}
    \end{minipage}
    &
    \begin{minipage}[c]{\x\textwidth}
       \centering 
        \subfloat[Periodic, $X_{bc}$ = 15 mm]{\includegraphics[trim={10cm 0 14cm 0},clip,width=\textwidth]{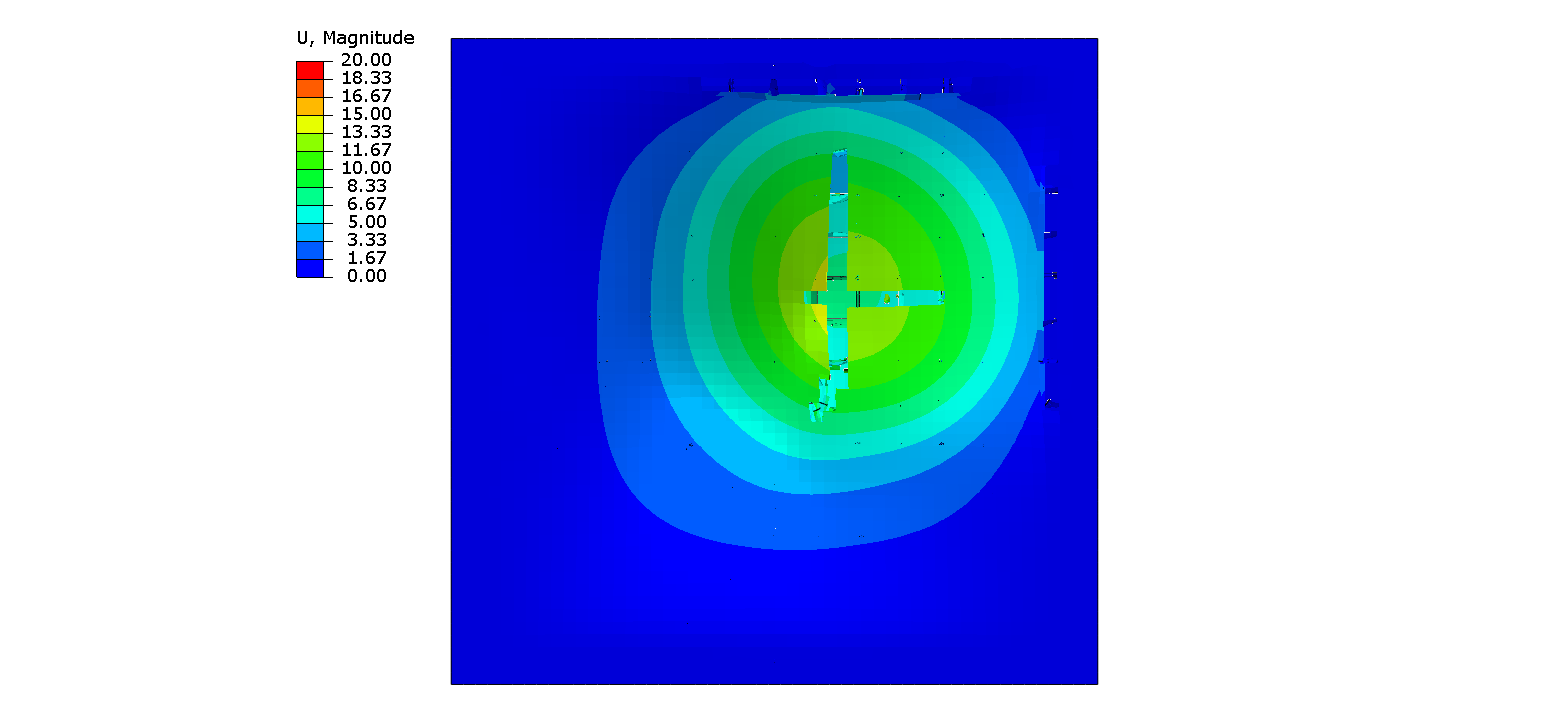}
        \label{fig:d9}}
    \end{minipage}

    \end{tabular}
    \caption{Top face sheet displacement of different lattice core designs with different blast center locations. All displacements color bars are identical. Simulations correspond to the original blast center location (15,15) mm are shown in the third column.}
    \label{shift_disp_comp}
\end{figure}

Finally, to evaluate the robustness of the periodic design, we gradually moved the location of the blast center from (15,15) mm to (-15,-15) mm, while maintaining the same blast load magnitude. We tested the initial design, and the final optimized designs from the free optimization (\sref{sec:ex2}) and the periodic optimization (\sref{sec:ex3}). The plot of the normalized energy quantities for the sandwich panel system (including the face sheets) is presented in \fref{shifted_loc}. Selected surface plots of displacement of the top face sheet are presented in \fref{shift_disp_comp}.

From \fref{shifted_loc}, we see that the periodic optimized design only performed better than the initial lattice design for the last two cases, where the blast center remained close to its original location. As the blast center shifted further away from its original position, the initial design outperformed the periodic optimized, showing higher stiffness (indicated by lower external work) and less plastic and damage dissipation. Interestingly, we noticed that the free optimized design, although optimized with a blast center location of (15,15) mm, remained largely better than both the initial and periodic optimized designs, especially in terms of the high structure stiffness and low material damage. The free optimized design was only inferior to the other designs for a blast center location of (-15,-15) mm, which is reasonable since that region of the face sheet was poorly supported by the lattice core (see the lattice wall thickness distribution in \fref{fig:d4_opt}). The surface displacement plots in \fref{shift_disp_comp} further confirms the trends seen in \fref{shifted_loc}. We observe that both optimized designs are effective in reducing the top face sheet displacements at the original blast center location (third column of \fref{shift_disp_comp}). However, as the blast shifted to the center of the face sheet (center column of \fref{shift_disp_comp}), the periodic optimized design became less effective and started showing top face sheet fracture, which was absent from the free optimized design. When the blast center shifted to (-15,-15) mm (first column in \fref{shift_disp_comp}), we notice a new failure mode emerged for the free optimized design, where both the top and bottom face sheets sheared off the fixed edges due to the lack of underlying lattice core support.

Evidence presented in both \fref{shifted_loc} and \fref{shift_disp_comp} indicate that the periodic optimized design is inferior to the free optimized design. This observation is reasonable for the specific setup of these two examples, where the periodic optimization case had a much smaller design space compared to the free optimization case (84 vs. 364 designable lattice walls). Optimized design performance typically deteriorates as the available design space become smaller and more restricted. If we were to increase the design space by increasing the number of lattice walls in one periodic unit cell, the optimized design is expected to provide better performance than the current one.

\section{Conclusions and future work}
\label{sec:conc}
In this work, a heuristic topology optimization framework for thin-walled lattice structure is proposed and tested. The framework relies on a heuristic thickness update scheme that is based on the idea of equalization of energy density across all lattice walls. The thickness update is achieved by solving a nonlinear optimization defined by the update scheme. Two novel thickness update schemes are presented, one being a direct minimization of the scatter in energy density, and the other is inspired by the BESO method. The framework allows the specification of minimum and maximum allowable wall thicknesses for manufacturing considerations, and lattice walls with thickness below the minimum threshold are removed from the structure. The framework also allows for the generation of periodic lattice structures through the use of periodic aggregation and scattering operators, which ensures the periodicity of the generated design. The proposed framework terminates either when the thickness change is below a minimum threshold, or when the maximum number of design iterations is reached.

Three numerical examples are presented to demonstrate the capabilities of the proposed framework. In the first example, we performed a benchmark test with the previous HCA-based topology optimization framework by \cite{hunkeler2014topology}, where both thickness update schemes were capable of generating designs with higher specific energy absorption than the optimized designs in \cite{hunkeler2014topology}. In the second example, we performed topology optimization for a sandwich panel under blast loading, considering material damage behavior. It was found that the proposed thickness update schemes were effective in reducing damage and increasing stiffness of the sandwich panel, and when the load magnitude is increased, increasing specific energy absorption. The last example follows from the previous one, except we generated a periodic lattice. We observed improved stiffness and reduced damage compared to the initial design, but the performance of the periodic design was inferior to the free design. In all three examples, the proposed framework was able to greatly improve key performance metrics such as specific energy absorption, stiffness and material damage in about 25 design iterations, and in the first example, it is more efficient than the HCA framework proposed in \cite{hunkeler2014topology} and the commercial software LS-OPT. Although all simulations presented in this work were conducted in Abaqus \cite{Abaqus2021}, the proposed framework can be easily adapted to work with other finite element software. 

We conclude that the proposed topology optimization framework can effectively generate optimized thin-walled lattice structures under different loading conditions. The high effectiveness and the ability to handle complex loading and material behavior renders our framework a suitable tool to generate optimized lattice structures for energy absorbing applications and for improving crashworthiness of components. 

In future work, we will focus on validating the generated designs through experiments. In addition, it is known in the literature that the material properties of additively manufactured materials are typically depend on the feature size (e.g., wall thickness and strut diameter) \cite{phutela2019effects,barba2020size,alghamdi2021buckling}. Therefore, it is of interest to leverage the current framework to study how the size dependence affects the optimized lattice designs at different length scales. The current framework also provides a physics-informed way to generate designs with high performance for training of machine learning frameworks, which can be a more suitable alternative to the random combinatorial generation method used in our previous work \cite{he2022exploring}.

\section*{Replication of results}
The data and source code that support the findings of this study can be found at: \url{https://github.com/Jasiuk-Research-Group/LatticeOPT}. \textcolor{red}{Note to editor and reviewers: the link above will be made public upon the publication of this manuscript. During the review period, the data and source code can be made available upon request to the corresponding author.}

\section*{Conflict of interest}
The authors declare that they have no conflict of interest.

\section*{Acknowledgements}
We (I. J.) acknowledges the support of the Army Research Office contract (No. W 911NF-18-2-0067) and the National Science Foundation grant (MOMS-1926353).

\section*{CRediT author contributions}
\textbf{Junyan He}: Conceptualization, Methodology, Software, Formal analysis, Investigation, Data Curation, Writing - Original Draft. \textbf{Shashank Kushwaha}: Methodology, Investigation, Writing - Original Draft. \textbf{Diab Abueidda}: Supervision, Writing - Review \& Editing. \textbf{Iwona Jasiuk}: Supervision, Resources, Writing - Review \& Editing, Funding Acquisition.

\bibliographystyle{unsrtnat}
\setlength{\bibsep}{0.0pt}
{\scriptsize \bibliography{References.bib} }
\end{document}